\documentclass[11pt]{amsart}

\usepackage{amsmath,amsthm, amscd, amssymb, amsfonts}
\usepackage[all]{xy}
\usepackage{enumitem}
\usepackage{array}  
\newcolumntype{L}{>{$}l<{$}} 
\newcolumntype{R}{>{$}r<{$}}
\newcolumntype{C}{>{$}c<{$}}
\usepackage{fontsmpl}
\usepackage{srcltx}
\usepackage[toc,page]{appendix}

\usepackage[T2A,T1]{fontenc}

\newcommand{\bsl}{\boldsymbol{\lambda}}

\allowdisplaybreaks

\newcommand{\supera}[2]{{\mathbf A}(#1|#2)}

\newcommand{\toba}{\mathcal{B}}
\newcommand{\wtoba}{\widetilde{\toba}}

\newcommand{\cL}{\mathcal{L}}

\newcommand{\cR}{\mathcal R}

\newcommand{\cI}{{\mathcal I}}
\newcommand{\cJ}{{\mathcal J}}
\newcommand{\cJt}{\widetilde{{\mathcal J}}}
\newcommand{\cG}{{\mathcal G}}
\newcommand{\cZ}{{\mathcal Z}}

\newcommand{\nombre}{{\mathtt{HV}_1}}
\newcommand{\runow}{V_{12}}
\newcommand{\rdosw}{V_{112}}
\newcommand{\FK}{\mathtt{FK}}

\usepackage{color}

\numberwithin{equation}{section}\theoremstyle{plain}

\newtheorem{step}{Step}
\newtheorem{theorem}{Theorem}[section]

\newtheorem{lemma}[theorem]{Lemma}

\newtheorem{pro}[theorem]{Proposition}
\newtheorem*{claimo}{Claim}
\newtheorem*{lemma*}{Lemma}

\newtheorem{claim}{Claim}

\theoremstyle{definition}

\newtheorem{definition}[theorem]{Definition}

\newtheorem{example}[theorem]{Example}

\theoremstyle{remark}
\newtheorem{rem}[theorem]{Remark}
\newtheorem{remark}[theorem]{Remark}

\newcommand{\ydh}{{}^{H}_{H}\mathcal{YD}}

\newcommand{\ydg}{{}^{\ku\Gamma}_{\ku\Gamma}\mathcal{YD}}

\newcommand\Ss{\mathcal{S}}

\newcommand\G{\mathbb{G}}

\newcommand\Gp{\mathbb{G}'}

\newcommand\Ib{\mathbb{I}}
\newcommand\I{\mathbb{I}}
\newcommand\id{\operatorname{id}}

\newcommand\Alg{\operatorname{Alg}}
\newcommand\Cleft{\operatorname{Cleft}}

\newcommand\GK{\operatorname{GK-dim}}

\newcommand\End{\operatorname{End}}

\newcommand\gr{\operatorname{gr}}
\newcommand\co{\operatorname{co}}

\newcommand\ad{\operatorname{ad}}

\newcommand{\supp}{\operatorname{supp}}

\newcommand\ku{\Bbbk}

\newcommand\bq{\mathfrak{q}}

\newcommand\ot{\otimes}

\newcommand\s{\mathbb{S}}

\newcommand\N{\mathbb{N}}

\newcommand\mT{\mathcal{T}}

\newcommand\mD{\mathcal{D}}

\newcommand\mP{\mathcal{P}}
\newcommand\mO{\mathcal{O}}

\newcommand\cA{\mathcal{A}}
\newcommand\mH{\mathcal{H}}

\newcommand{\att}{\mathtt{a}}
\newcommand{\ytt}{\mathtt{y}}

\newcommand{\J}{{\mathcal J}}

\newcommand{\cE}{{\mathcal E}}
\newcommand{\z}[1]{x_{#1 4}}
\newcommand{\hu}[1]{x_{#1 4}}

\newcommand{\wx}{\widetilde{x}}
\newcommand{\wu}[1]{\widetilde{x}_{#1 4}}
\newcommand{\wz}[1]{\widetilde{x}_{#1 4}}
\newcommand{\rel}{\mathtt{r}}
\newcommand{\relp}{\mathtt{p}}

\newcommand\pf{\begin{proof}}
\newcommand\epf{\end{proof}}

\begin{document}

\title[Pointed Hopf algebras over non abelian groups, I]
{Pointed Hopf algebras over non abelian groups with decomposable braidings, I}

\author[Iv\'an Angiono and Guillermo Sanmarco]
{Iv\'an Angiono and Guillermo Sanmarco}

\thanks{ S. was partially supported by CONICET, Secyt (UNC), the MathAmSud project GR2HOPF}

\address{ Facultad de Matem\'atica, Astronom\'{\i}a y F\'{\i}sica,
Universidad Nacional de C\'ordoba. CIEM -- CONICET. %\newline 
Medina Allende s/n (5000) Ciudad Universitaria, C\'ordoba,
Argentina}
\email{(angiono|gsanmarco)@famaf.unc.edu.ar}

\begin{abstract}
We describe all finite-dimensional pointed Hopf algebras whose infinitesimal braiding is a fixed Yetter-Drinfeld module decomposed as the sum of two simple objects: a point and the one of transpositions of the symmetric group in three letters. We give a presentation by generators and relations of the corresponding Nichols algebra and show that Andruskiewitsch-Schneider Conjecture holds for this kind of pointed Hopf algebras.
\end{abstract} 
\maketitle

\section{Introduction}

The theory of Hopf algebras has grown recently, 
focused especially in the finite dimensional case. The idea of classifying these algebras 
leads to separate into families that share certain invariants, mainly associated 
with their coalgebra structure. The first invariant to take in care is the coradical (the sum of all simple subcoalgebras).
The first distinguished family is that of cosemisimple Hopf algebras, those which coincide with the coradical. A second important family is constituted by pointed Hopf algebras
(those whose coradical coincides with the subalgebra generated by group-like elements) and was intensively studied after the introduction of quantized enveloping algebras by Drinfeld and Jimbo.

A turning point in the classification of pointed Hopf algebras was given by the introduction of the Lifting Method by Andruskiewitsch and Schneider \cite{AS-adv}: a sequence of steps proposed to systematically describe all the deformations of a graded Hopf algebra attached to a fixed finite group $\Gamma$ (giving the coradical) and a Yetter-Drinfeld module $V$ over $\Gamma$ (describing the \emph{infinitesimal braiding}). At this point an algebra $\toba(V)$ universally contrusted from $V$ appears: the so-called \emph{Nichols algebra} of $V$, cf. \S \ref{subsec:Nichols}.
The steps of this method are the following:
\begin{enumerate}[leftmargin=*,label=\rm{(\alph*)}]
\item\label{item:liftingstep1} Classify all Yetter-Drinfeld modules $V$ over $\Gamma$ such that $\dim\toba(V)<\infty$.
\item\label{item:liftingstep2} Give a presentation by generators and relations of $\toba(V)$.
\item\label{item:liftingstep3} Check if, given $H$ with $H_0\simeq \Bbbk\Gamma$, the associated graded Hopf algebra satisfies $\gr H\simeq \toba(V)\# \Bbbk\Gamma$ (generation in degree one problem).
\item\label{item:liftingstep4} Compute all Hopf algebras obtained as deformations of $\toba(V)\# \Bbbk\Gamma$.
\end{enumerate}
The Lifting Method was succefully applied to classify finite-dimensional pointed Hopf algebras over abelian groups $\Gamma$ of order not divisible by small primes \cite{AS-Annals}.
After the appeareance of this result, Masouka \cite{Ma} was able to prove that all these pointed Hopf algebras are cocycle deformations of their associated graded Hopf algebras.

For general finite abelian groups $\Gamma$, Heckenberger \cite{H-rank2root} gave a complete answer of \ref{item:liftingstep1} and Angiono \cite{An-diagonal} finished \ref{item:liftingstep2} and \ref{item:liftingstep3}. The final step was recently completed in \cite{AG}, based on a strategy to construct Hopf algebras using cocycle deformations developed in \cite{A+,AAG}.

If we turn towards non-abelian groups, the first step \ref{item:liftingstep1} was almost completed \cite{HV-rank2,HV-rank>2}. The authors gave all non-simple Yetter-Drinfeld modules over non-abelian groups such that the associated Nichols algebras are finite-dimensional. They describe their root systems and compute the dimensions. At the moment there is no general answers for the remaining steps of the Lifting Method. Anyway the method was applied successfully on certain pairs of non-abelian groups and infinitesimal braidings whose Nichols algebras are finite dimensional, see e.g. \cite{GV,GV2}

The present work starts the computation of all liftings of Nichols algebras of Yetter-Drinfeld modules classified in \cite{HV-rank2,HV-rank>2}. We fix a brainding (that we call $\nombre$) decomposable as the sum of two simple Yetter-Drinfeld modules, one related with the so-called Fomin-Kirillov algebra \cite{MS,AG} over the symmetric group in three letters and another of dimension one. We complete the answer for steps \ref{item:liftingstep2}, \ref{item:liftingstep3} and \ref{item:liftingstep4}. Indeed, we give a presentation by generators and relations of $\toba(\nombre)$ in \S \ref{sec:D3-1}. Next we give a positive answer to the generation-in-degree-one problem in \S \ref{sec:preNichols}: the unique finite-dimensional post-Nichols algebra of $\nombre$ is the Nichols algebra $\toba(\nombre)$ itself. We also introduce a \emph{distinguished} pre-Nichols algebra of $\nombre$, whose behavior is analogous to those pre-Nichols algebras in \cite{An-distinguished} for braidings of diagonal type. Finally, the main result of this paper is contained in \S \ref{sec:liftings} and gives a complete list of all finite-dimensional pointed Hopf algebras whose group of group-like elements is a group $\Gamma$ such that $\nombre$ admits a \emph{principal realization} over $\Gamma$. We follow the strategy of \cite{A+,AAG} and prove at the same time that all of them are cocycle deformations of $\toba(\nombre)\# \Bbbk\Gamma$.

\subsection*{Acknowledgements}

We thank Nicol\'as Andruskiewitsch for posing us the problem and contributing with ideas for its solution. Many interesting discussions with him helped us to improve the paper.

\section{Preliminaries}
\subsection{Notation}
If $\ell < \theta \in\N_0$, then we set $\I_{\ell, \theta}=\{\ell, \ell +1,\dots,\theta\}$, $\I_\theta 
= \I_{1, \theta}$.

Throughout this work $\Bbbk$ denotes an algebraically closed field of characteristic zero. All  vector spaces, algebras and tensor products  are over $\ku$.
Let $\G_N$ be the group of roots of unity of order $N$ in $\ku$ and $\G_N'$ the subset of primitive roots of order $N$;
$\G_{\infty} = \bigcup_{N\in \N} \G_N$ and $\G'_{\infty} = \G_{\infty} - \{1\}$.

Let $A$ be an algebra. The ideal generated by a subset $(a_i)_{i\in I}$ is denoted by $\langle a_i: i\in I\rangle$  while the subalgebra generated by this set is denoted by $\Bbbk\langle a_i: i\in I\rangle$.

Let $H$ be a Hopf algebra (or a Hopf algebra in a braided tensor category). We will use the Sweedler notation for the coalgebra structure and for comodules; e.g. the coproduct of $h\in H$ is denoted $\Delta(h)= h_{(1)}\ot h_{(2)}$. Let $G(H)=\{x\in H^{\times}: \Delta(x)=x\ot x\}$ be the group of group-like elements. For each $g,h\in G(H)$, $\mP_{g,h}(H)$ is the subspace of $(g,h)$-primitive elements; in particular $\mP(H):= \mP_{1,1}(H)$ is the set of primitive elements.

\subsection{Braided vector spaces, Yetter-Drinfeld modules and braided Hopf algebras}
A pair $(V,c)$, where $V$ is a vector space and $c\in \operatorname{GL}(V\ot V)$, is a braided vector space if $c$ satisfies the braid equation:
\begin{align*}
(c\ot \id)(\id\ot c)(c\ot \id)=(\id\ot c)(c\ot \id)(\id\ot c).
\end{align*}

\medbreak

Let $H$ be a Hopf algebra with bijective antipode. Attached to $H$ there is a braided tensor category $\ydh$. An object of $\ydh$, called
a \emph{Yetter-Drinfeld module over $H$}, is a left $H$-module and left $H$-comodule $V$ such that
\begin{align}\label{eq:yd-comp}
\delta (h\cdot v)= h_{(1)}v_{(-1)}\Ss(h_{(3)}) \ot h_{(2)} \cdot v_{(0)} && \text {for all }  h\in H, & v\in V.
\end{align}
The monoidal structure on $\ydh$ comes from the usual tensor product of modules and comodules over a Hopf algebra. 
The left dual of a finite dimensional $V \in \ydh$ is the vector space $V^*$, where $H$ acts via the antipode $(h\cdot f)(v)=f(\Ss(h)\cdot v)$
and the coaction $\delta(f)= f_{(-1)} \ot f_{(0)}$ is determined by 
\begin{align*}
f_{(-1)}  f_{(0)} (v) &= \Ss^{-1} (v _{(-1)})  f (v_{(0)}), & f&\in V^*, \quad v \in V.
\end{align*}
For $V, U\in \ydh$, the braiding $c_{V, U} \colon V\ot U \to U\ot V$ is given by
\begin{align}\label{eq:yd-braid}
c_{V,U}(v\ot u) &= v_{(-1)}\cdot u \ot v_{(0)}, & v\in V, \, \, & u\in U.
\end{align}
The pair $(V, c_{V,V})$ is a braided vector space. 

Reciprocally, let $(V,c)$ be a braided vector space. A \emph{realization} of $V$ over $H$ is an structure on $V$ of Yetter-Drinfeld module over $H$ such that the braiding $c$ is the categorical braiding $c_{V,V}$. The realization is \emph{principal} if there exists a basis $\{v_i\}_{i\in I}$ of $V$ and $\{g_i\}_{i\in I}\in G(H)$ such that $\delta(v_i)=g_i\ot v_i$, $i\in I$.

\medbreak

There is a monoidal structure in the category of algebras in $\ydh$. 
Namely, if $(A, \mu_A)$ and $(B, \mu_B)$ are algebras in $\ydh$, then $A\underline{\ot}B=(A\ot B, \mu_{A\ot B})$ also is, where the multiplication is given by
\begin{align}
\mu_{A\ot B}=(\mu_A\ot \mu_B)(\id_A\ot c_{A,B} \ot \id_B).
\end{align}

\bigbreak

We are mainly interested in Yetter-Drinfeld modules over a group algebra $\Bbbk\Gamma$.
Let $V$ be simultaneously a left $\Bbbk \Gamma$-comodule and a left $\Bbbk\Gamma$-module.
Then $V=\oplus_{g\in \Gamma}V_g$, where $V_g=\{v\in V \colon \delta(v)=v\ot g\}$. 
In this setting the Yetter-Drinfeld compatibility \eqref{eq:yd-comp} is read as
\begin{align}\label{eq:ydg-braid}
\delta (g\cdot v)= ghg^{-1} \ot g \cdot v && \text { for all }  g, h\in \Gamma, & v\in V_h .
\end{align}
Thus $V \in \ydg$ means that $g\cdot V_h \subset V_{ghg^{-1}}$ for all $g, h \in \Gamma$.

We recall that the \emph{support} of $V$ is the subset 
\begin{align}\label{eq:supp-YD}
\supp V := \{ g\in\Gamma: V_g\neq 0 \}  \subset \Gamma.
\end{align}
We denote by $\widehat{\Gamma}$ the group of characters of $\Gamma$.
For each $\chi\in\widehat{\Gamma}$,
\begin{align}\label{eq:supp-YD}
V^{\chi} := \{v\in V: g\cdot v= \chi(g)v \text{ for all }g\in\Gamma \}.
\end{align}
We also use the following notation: $V_g^{\chi}=V_g \cap V^{\chi}$.

We recall that a  Hopf algebra in $\ydh$ is a collection $(R, \mu, \Delta, \Ss)$, where $R$ is an object $\ydh$ with structures $(R, \mu)$ of algebra 
in $\ydh$ and $(R, \Delta)$ of coalgebra in $\ydh$ that are compatible in the sense that 
$\Delta \colon R \to R\underline{\ot}R$ and $\varepsilon \colon R \to \ku$ are algebra maps, 
and such that $\Ss$ is a convolution inverse of the identity of $R$.
In this case, the left adjoint action of $R$ on itself is the linear map $\ad_c \colon R \to \End(R)$,
\begin{align*}
(\ad_cx)y&=\mu(\mu\ot \Ss)(id\ot c)(\Delta \ot id)(x\ot y), &  x, y \in R,
\end{align*}
so the action of a primitive element $x \in R$ is
\begin{align*}
(\ad_cx)y&=xy - (x_{(-1)}\cdot y)x_{(0)}, &   y \in R.
\end{align*}

\subsection{Racks}\label{subsec:racks-realizations}

We recall now a prototypical example of Yetter-Drinfeld modules over groups. A \emph{rack} $X=(X,\rhd)$ is a pair where $X$ is a non-empty set and  $\rhd:X\times X\rightarrow X$ is an operation such that
\begin{align*}
x\rhd (y\rhd z) &=(x\rhd y)\rhd (x\rhd z), & \text{for all }&x,y,z\in X
\end{align*}
and the maps $\phi_x\colon X\longmapsto X$, $\phi_x(y)=x\rhd y$, $y\in X$,
are bijective for all $x\in X$. A \emph{quandle} is a rack such that $x\rhd x=x$ for all $x\in X$.

Let $\Gamma$ be a group. An example of a rack (which is moreover a quandle) is given by $X$ a union of conjugacy classes inside $\Gamma$ and the operation given by $g \rhd h=ghg^{-1}$, $g,h\in \mO$. 

A \emph{2-cocycle} on $X$ is a function  $\bq:X\times X\rightarrow\Bbbk^\times$, $(x,y)\mapsto q_{x,y}$, such that
\begin{align}\label{eq:cocycle-rack-condition}
q_{x,y\rhd z}q_{y,z}&=q_{x\rhd y,x\rhd z}q_{x,z}, & \text{for all }&x,y,z\in X.
\end{align}
Let $(X,\rhd)$ be a rack and $\bq$ a 2-cocycle on $X$.
Let $V$ be the vector space with basis $\{v_x\colon x\in X\}$. Let $c: V\ot V\to V\ot V$ be the linear map
\begin{align}\label{eqn:rack-braiding}
c(v_x\ot v_y)&=q_{x,y}v_{x\rhd y}\ot v_x, & x,y&\in X.
\end{align} 
The pair $(V,c)$ is a braided vector space denoted by $V(X,\bq)$. 

\medspace

Next we recall the definition of a \emph{principal realization} of a braided vector space $V(X,\bq)$ (here we restrict to the case of group algebras as in \cite[Definition 3.2]{AG}, see \cite[\S 4]{GV2} for the case of an arbitrary Hopf algebra). It is given by the data $(\cdot, (g_i)_{i\in X}, (\chi_i)_{i\in X})$ consisting of:
\begin{itemize}[leftmargin=*]
\item an action $\cdot:G\times X\to X$ of $G$ on $X$; 
\item a family of elements $g_i \in G$, $i\in X$,  such that $g_{h\cdot i} = hg_ih^{-1}$ and $g_i\cdot j=i\rhd j$, for all $i, j \in X$, $h\in G$; 
\item a family of functions $\chi_i:G\to\Bbbk^\times$, $i\in X$, called $1$-cocycles, such that $\chi_i(g(j))=q_{ji}$, $i,j\in X$, $\chi_i(ht)=\chi_i(t)\chi_{t\cdot i}(h)$, for all $i\in X$, $h,t\in G$.
\end{itemize}
This data defines a Yetter-Drinfeld structure on $V(X,q)$ by 
\begin{align}
\delta(v_i) &= g_i\ot v_i, & g\cdot v_i&=\chi_i(g)v_{g\cdot i}, & 
i\in X,& \,\, g\in G.
\end{align}

\subsection{Nichols algebras}\label{subsec:Nichols}
Let $V \in \ydh$. The Nichols algebra of $V$ is the unique quotient of the tensor algebra $\toba(V)= T(V)/ \cJ (V)$  in $\ydh$ that is a graded 
braided Hopf algebra  and such that $V$ is the space of primitive elements. 
Among the several descriptions of the universally defined ideal $\cJ (V) = \oplus_{n\geq2} \cJ_n(V)$, we recall the most useful for our purposes: $\cJ(V)$ is the unique maximal element in the class of graded Hopf ideals and Yetter-Drinfeld submodules of $T(V)$ with zero components on degree $0$ and $1$.

\bigbreak

We recall that a \emph{pre-Nichols algebra} $\wtoba$ of $V$ is 
a graded braided Hopf algebra $\wtoba=\oplus_{n\geq0} \wtoba^n$ such that $\wtoba^0=\Bbbk1$, $\wtoba^1=V$ and $\wtoba$ is generated as an algebra by $V$; that is, a graded braided Hopf algebra, intermediate between $T(V)$ and $\toba(V)$. By definition, the inclusion $V\hookrightarrow \wtoba$ induces graded Hopf algebra epimorphisms $T(V)\twoheadrightarrow \wtoba$ and $\wtoba\twoheadrightarrow \toba(V)$. 

Fix $(x_i)_{i\in I}$ a basis of $V$, then we use the following notation:
\begin{align}\label{eq:adjunta-iterada}
x_{i_1 \dots i_n} &:= (\ad_c x_{i_1}) \dots (\ad_c x_{i_{n-1}}) x_{i_n}, &
n\ge 2, & \, i_1, \dots, i_n \in I.
\end{align}

Dually, a \emph{post-Nichols algebra} $\cL$ of $V$ is a graded braided Hopf algebra, intermediate between $T^c(V)$, the free associative coalgebra of $V$, and $\toba(V)$; that is, a graded braided Hopf algebra $\cL=\oplus_{n\geq0} \cL^n$ such that $\cL^0=\Bbbk1$, $\cL^1=V$ and $\mathcal{P}(\cL)=\cL^1$.

\bigbreak

Suppose $V$ has a basis $(x_i)_{ i \in \I}$ where the coaction is given by $\delta(x_i)=g_i \ot x_i$ for some $(g_i)_{ i \in \I} \subset G(H)$. 
Let $(f_i)_{ i \in \I}$ be the dual basis of $(x_i)_{ i \in \I}$. Define a \emph{skew derivation} $\partial_i \colon T(V) \to T(V)$ by 
$\partial_i (1) = 0$,   $\partial_i |_V = f_i$ and recursively
\begin{align*}
\partial_i (xy) &= x  \partial_i (y) + \partial_i (x) g_i \cdot y, & x, y &\in T(V).
\end{align*}
It descends to every pre-Nichols algebra $\toba$ and gives a criterion to characterize Nichols algebras:
\begin{align*}
\toba&=\toba(V) && \iff &\bigcap_{i\in \I} \ker \partial_i  & = \Bbbk 1.
\end{align*}
In particular every $x\in \toba$ such that $\partial_i(x)=0$ for all $i\in\I$ projects to 0 onto the Nichols algebra $\toba(V)$.

\subsection{Cocycle deformations, cleft objects}
A \emph{normalized $2$-cocycle} on a Hopf algebra $H$ is a convolution invertible map $\sigma \colon H \ot H \to \ku$ such that
\begin{align*}
\begin{aligned}
\sigma (x,1)= \sigma (1,x)&=\varepsilon (x), \\
\sigma( x_{(1)}, y_{(1)}) \sigma( x_{(2)} y_{(2)}, z) &= \sigma( y_{(1)}, z_{(1)}) \sigma( x, y_{(2)} z_{(2)}),
\end{aligned}
&&  x, y, z \in H.
\end{align*}
This allows to modify the multiplication of $H$  to a new associative multiplication $ \cdot_{\sigma} \colon H \ot H \to \ku$ given by
\begin{align*}
 x \cdot_{\sigma} y &= \sigma( x_{(1)}, y_{(1)})  x_{(2)} y_{(2)}  \sigma ^ {-1}( x_{(3)}, y_{(3)}), & x, y \in H.
\end{align*}
The antipode $\Ss$ of $H$ can also be modified via $\sigma$ to produce certain $\Ss_{\sigma}$ in such a way 
that $H_{\sigma} := (H, \cdot_{\sigma}, 1, \Delta, \varepsilon, \Ss_{\sigma})$ is a Hopf algebra.
We denote the set of $2-$cocycles on $H$ by $Z^2(H, \ku)$.

A  \emph{cocycle deformation} of $H$ is a Hopf algebra isomorphic to some $H_{\sigma} $.
\medbreak
A \emph{right cleft object} for $H$ is a right $H$-comodule algebra $C$ that satisfies $C^{\co H} = \ku$ and admits a convolution invertible 
comodule map $\gamma \colon H \to C$. In this case, $\gamma$ can be normalized to $\gamma(1) = 1$ and is called a \emph{section}.
We denote by $\Cleft (H)$ the set of iso-classes of cleft objects.
\medspace

There exists an analogous notion of left cleft object. Also, given $K,H$ two Hopf algebras, $A$ is a $(K,H)$-bicleft object if $A$ is simultaneously a right $H$-cleft object and a left $K$-object, and the two coactions commute.

\medbreak
These two notions are equivalent. Indeed each cocycle $\sigma$ gives raise to a right cleft object $A$ such that $A$ is a  $(H_{\sigma},H)$bicleft object. Reciprocally for each right $H$-cleft object $A$ there exists a Hopf algebra $L:=L(A,H)$, called the
Schauenburg algebra, such that $A$ is $(L,H)$-bicleft and $L(A,H) \simeq H_{\sigma}$ for some $\sigma$, see \cite{S} for details; $L$ is unique up to isomorphism.

\subsection{Liftings of the Fomin-Kirillov algebra $\FK_3$}\label{subsec:FK}
Let $\mathcal{O}_2^3$ be the rack of transpositions in $\s_3$.
Let $(V,c)$ be the braided  vector space determined by $\mathcal{O}_2^3$ with cocycle $\bq \equiv -1$.
Thus $V$ has a basis $(x_i)_{i\in \Ib_3}$ such that
\begin{align}
c(x_i \ot x_j) &= -x_{2i-j} \ot x_i,& i,j \in & \I_3,
\end{align}
where $2i-j$ means the class modulo 3.

By \cite{FK,MS}, the Nichols algebra ${\FK}_{3} := \toba(V)$ is presented by generators $x_{i}$, $i \in \Ib_3$, with relations
\begin{align}\label{eq:rels-FK3}
\begin{aligned}
x_{i}^{2}&=0,&  & i  \in \Ib_3,\\
x_{1}x_{2} + x_{3}x_{1} +x_{2}x_{3} &=0,&  
&x_{2}x_{1} + x_{1}x_{3} +x_{3}x_{2} =0.
\end{aligned}
\end{align}
Also, $\dim \FK_{3}=12$. Indeed, the following set is a basis of  $\FK_{3}$:
\begin{align*}
\{x_1^{n_1} x_2^{n_2} x_3^{n_3}\colon 0\leq n_1, n_2, n_3\leq 1 \} \cup  
\{x_1^{n_1} x_2x_1 x_3^{n_3}\colon 0\leq n_1, n_3\leq 1 \}.
\end{align*}

\medspace

Let $\Gamma$ be a group such that $V$ admits a principal realization $(g_i, \chi_i)_{i\in \I_3}$ over $H=\Bbbk\Gamma$, see \S \ref{subsec:racks-realizations}.
% where $g_i\in G(H)$ are such that $x_i\in V_{g_i}$. 
For example we may take $\Gamma=\s_3$. By \cite[Theorem 1.2]{GV2} $L$ is a lifting of $\FK_3$ over $H$ if and only if $L$ is (isomorphic to) the quotient of $T(V)\# H$ by the relations
\begin{align}\label{eq:rels-FK3-liftings}
x_{1}^{2}&=\lambda_1(1-g_i^2), &
x_{1}x_{2} + x_{3}x_{1} +x_{2}x_{3} &=\lambda_{2}(1-g_1g_2),
\end{align}
where the \emph{deformation parameters} $\lambda_1, \lambda_2 \in \ku$ satisfy the constraints
\begin{align}\label{eq:lambda-conditions-FK3}
\begin{aligned}
\lambda_1 &= 0 \text { if either } \chi_i^2 \neq \varepsilon \text { or } g_i^2 = 1 \text { for some } i\in\I_3, 
\\
\lambda_2 &= 0 \text { if either } \chi_i \chi_j  \neq \varepsilon \text { or } g_i g_j = 1 \text { for some } i\neq j\in\I_3.
\end{aligned}
\end{align}
Notice that \eqref{eq:rels-FK3} implies that the following relations also hold:
\begin{align*}
x_{i}^{2}&=\lambda_1(1-g_i^2), \,  i  \in \Ib_{2,3}, &
x_{2}x_{1} + x_{1}x_{3} +x_{3}x_{2} &=\lambda_{2}(1-g_2g_1).
\end{align*}

\section{The Nichols algebra $\toba(\nombre)$}\label{sec:D3-1}

In this section we introduce the Nichols algebra $\toba(\nombre)$. We recall the definition of principal realization for the associated braided vector space, describe the root system of $\nombre$ and give a (minimal) presentation by generators and relations of this Nichols algebra.

\subsection{The problem}
Our aim is to study the Nichols algebra $\toba(\nombre)$, cf. \cite[Example 1.10]{HV-rank2}.
Here $\nombre=V_1 \oplus V_2$, where $V_1$ is as in \S \ref{subsec:FK} and $V_2$ has dimension one: we fix $x_4\in V-0$. 
Let $q_1,q_2 \in \ku^{\times}$ be such that 
\begin{align}\label{eq:omega}
\omega := -q_1q_2 \in \Gp_3.
\end{align}  
The braiding of $\nombre$ is determined by
\begin{itemize}[leftmargin=*]
\item $(V_1,c)$ is a braided vector subspace;

\item $V_2$ is a point with label $-\omega^2$, so $c(x_4\ot x_4)=-\omega^2 x_4\ot x_4$;

\item the braiding between them is
\begin{align}\label{eq:D3-1}
c(x_i \ot x_4) &= q_1  x_4 \ot x_i, & 
c(x_4 \ot x_i) &= q_2 x_i \ot x_4, & i&\in \I_3.
\end{align}
\end{itemize}

The braided vector space $\nombre$ is associated to the rack $\mathcal{O}_2^3 \times \{4\}$ and the matrix 
$\mathfrak q:=\begin{pmatrix} -1 & q_1 \\ q_2 & -\omega^2\end{pmatrix}$, cf. \cite[Example 38]{A-leyva}.

\medspace

The root system of $\nombre$ is of standard type $B_2$, cf. \cite{HS-adv,HV-rank2}. We fix the following notation for the set of positive roots: $\varDelta_+ = \{ 1,112,12,2 \}$

\subsection{Principal realizations} \label{subsec:Gamma3}

Let $\Gamma$ be a group. We assume that there exist $g_i\in \Gamma$, $i\in\I_4$, and $1$-cocycles $\chi_i:\Gamma\to\ku$, such that
\begin{align}\label{eq:Gamma3-rels}
\begin{aligned}
g_i^{\Gamma} &=\{g_1 , g_2, g_3 \}, &
g_4^{\Gamma}&=\{g_4 \} &
g_i^2 &=g_j^2, &
g_i g_j &= g_{2 i - j} g_i, &
i, &  j \in \I_3;
\end{aligned}
\\ \label{eq:cocycles}
\begin{aligned}
\chi _ i (g_j) &= \begin{cases} -1, & j \in \I_3, \\q_2,& j=4 \end{cases}, &  i& \in \I_3, &
\chi _ 4 (g_j) &= \begin{cases} q_1, & j \in \I_3, \\-\omega ^ 2, & j=4 \end{cases}.
\end{aligned}
\end{align}
Thus the conjugation induces an action $\triangleright \colon \Gamma \times \I_4 \to \I_4$ given by $g_{g \triangleright i}=g g_i g^{-1}$. Following \cite[Definition 3.2]{AG} we obtain a \emph{principal YD-realization} of $(\mathcal{O}_2^3 \times {4}, \mathfrak q )$ over $\Gamma$: $\nombre$ is realized in $^{\ku\Gamma}_{\ku\Gamma}\mathcal{YD}$ by defining
\begin{align}\label{eq:table-action}
\gamma \cdot x_i &= \chi_i(\gamma) x_{\gamma \triangleright i}, &
\delta(x_i) &= g_i \ot x_i,  & i\in & \mathbb{I}_4. 
\end{align}

\begin{example}
Let $\Gamma_3$ be the group generated by $\nu, \gamma, \zeta$ with relations
\begin{align*}
\gamma\nu &= \nu^2 \gamma, & \gamma \zeta & =\zeta\gamma,&
\nu\zeta &= \zeta\nu,& \nu^3&=1.
\end{align*}
It is isomorphic to the enveloping group of the quandle $\mathcal{O}_2^3$ \cite[Section 2]{HV-ii}. 
The elements $g_i := \gamma\nu^{i-1}$, $i\in \I_3$, and $g_4=\zeta$ generate $\Gamma_3$ and satisfy \eqref{eq:Gamma3-rels}.
There exist 1-cocycles $\chi_i \colon \Gamma_3 \to \ku ^{\times}$, $ i \in \mathcal{O}_2^3 \times {4}$, determined by
\begin{align}
\chi _ i (g_j) = \begin{cases} -1, & j \in \mathcal{O}_2^3 \\q_2,& j=4 \end{cases} \text { for } i \in \mathcal{O}_2^3,
\text{ and  }\chi _ 4 (g_j) = \begin{cases} q_1, & j \in \mathcal{O}_2^3 \\-\omega ^ 2, & j=4 \end{cases}.
\end{align}
Hence $\Gamma_3$ provides a principal YD-realization  of our quandle.

%We may also consider non-abelian groups $G$ with a surjective morphism $\pi \colon \Gamma_3 \twoheadrightarrow G$. Set $g_i := \pi (g_i)$, $i\in \I_4$. Since $G$ is non-abelian, $(g_i)_{i\in \I_4}$ are pairwise different.
\end{example}

We restrict now to a group $\Gamma$ as above and consider $\nombre \in \ydg$. 
For a shorter notation, we set
\begin{align*}
\runow &= (\ad_c V_1)V_2, &\rdosw &= (\ad_c V_1)^2 V_2.
\end{align*}

\begin{theorem}\label{thm:D3-1} \cite[Theorem 8.2]{HV-rank2}
The multiplication map 
\begin{align} \label{pbwD3-1}
\toba(V_2) \otimes \toba(\runow) \otimes \toba(\rdosw)\otimes \toba(V_1) \to \toba(\nombre)
\end{align}
is an isomorphism of $\N_0^2$-graded objects in $\ydg$.
Hence 
\begin{align*}
\dim \toba(\nombre) = 10368 = 3^42^7.
\end{align*}
\end{theorem}
\medbreak
\begin{remark}\label{rem:inclusion}
The algebra $\toba(V_2)$ is presented by a single generator $x_4$ and the relation $x_4^6=0$.
Since $V_1$ and $V_2$ are Yetter-Drinfeld submodules of $\nombre$, we have inclusions of Hopf algebras in $\ydg$
\begin{align*}
\toba(V_1) \hookrightarrow \toba(\nombre) \hookleftarrow \toba(V_2),
\end{align*}
thus the relations \eqref{eq:rels-FK3} and $x_4^6=0$ hold in $\toba(\nombre)$. 
\end{remark}

\subsection {Structure of $\runow = (\ad_c V_1)V_2$}  \label{subsec:runow}
Let
\begin{align}\label{eq:zi}
\z{i} &= (\ad_cx_i)x_4 \in \runow,  &  i&\in \Ib_3.
\end{align}
Some general features of the Yetter-Drinfeld module $\runow$ are given in \cite[Lemma 6.1]{HV-rank2}. Next we describe explicitly this structure for our realization.

\begin{lemma}\label{lema:zi} \begin{enumerate}[leftmargin=*,label=\rm{(\alph*)}]
\item\label{item:zi-props} The subspace $\runow$ is a Yetter-Drinfeld submodule of $\toba(\nombre)$ with basis $(\z{i})_{i\in \I_3}$. The action, coaction and braiding satisfy
\begin{align}\label{eq:zi-trenza}
&\begin{aligned}
g_i \cdot \z{j} &= - q_1 \z{2i-j \,}, & 
g_4 \cdot \z{j} &=- q_2\omega^2 \z{j}, & 
\\
\delta (\z{j})& = g_jg_4 \ot \z{j},& c(\z{i} \ot \z{j}) &= -\z{2i-j \,} \ot \z{i},
\end{aligned}&  i,j&\in \Ib_3.
\end{align}

\item\label{item:zi-map} The inclusion induces an algebra map $\toba(\runow) \to \toba(\nombre)$.
\end{enumerate}
\end{lemma}

\pf \ref{item:zi-props} Since $g_4$ in central in $G$ and $g_i \cdot x_4=q_1x_4$ for $i\in \Ib_3$, we have that
\begin{align*}
\z{i} = (\ad_cx_i)x_4 = x_ix_4 - q_1  x_4 x_i
\end{align*}
is homogeneous of degree $g_ig_4$. The set $(\z{i})_{i\in \I_3}$ spans $\runow$, hence it is enough to show that they are linearly independent. This is achieved computing their skew derivations:
\begin{align}\label{eq:zi-skew}
\partial_j(\z{i}) &= \delta_{ij} (g_j \cdot x_4 - q_1  x_4)=0, & 
\partial_4(\z{i}) &= x_i - q_1 g_4 \cdot x_i=-\omega^2x_i,
\end{align}
for  $i, j\in \Ib_3$, so they are linearly independent. 
The action of $\Gamma$ is given by
\begin{align*}
g_i \cdot \z{j} &=  g_i \cdot (x_jx_4 - q_1  x_4 x_j)=-q_1x_{2i - j}x_4+q_1^2x_4x_{2i - j}=- q_1 \z{2i-j \,}, \\
g_4 \cdot \z{j} &= -q_2\omega^2x_jx_4 + q_1 q_2\omega^2 x_4 x_j=-q_2\omega^2\z{j},
\end{align*}
and the formula for the braiding of $\runow$ follows from the action of $\Gamma$.

Now \ref{item:zi-map} follows by \cite[Theorem 2.6, Corollary 2.7(2)]{HS-JAlg}.
%By \ref{item:zi-props}, $\toba(\runow) \simeq\FK_3$.
%We show using skew derivations that $(\z{i})_{i\in \I_3}$ satisfies the relations 
%\eqref{eq:rels-FK3} in $\toba(\nombre)$.
%First, we show $\z{i}^2=0$. For $j \in \Ib_3$, it is clear $\partial_j(\z{i}^2)=0$, also $\partial_4(\z{i}^2)=$
%\begin{align*}
% = q_2\omega^4 x_i\z{i} - \omega^2 \z{i} x_i=q_2\omega x_i^2x_4+\omega^2 x_i x_4 x_i - \omega^2 x_i x_4 x_i +q_1\omega^2  x_4x_i^2 = 0.
%\end{align*}
%Next, clearly $\partial_j(\z{1}\z{2}+\z{3}\z{1}+\z{2}\z{3})=0$ for $j \in \Ib_3$. Now, if $i \neq j\in \I_3$, then
%\begin{align*}
%\partial_4(\z{i}\z{j})= q_2\omega^4 x_i\z{j} - \omega^2 \z{i} x_j = q_2 \omega  x_ix_jx_4 +q_1\omega^2 x_4x_ix_j .
%\end{align*}
%Hence  $\partial_4(\z{1}\z{2}+\z{3}\z{1}+\z{2}\z{3}) =$
%\begin{align*}
%&= -q_1 q_2^2 (x_1x_2+x_3x_1+x_2x_3)x_4 +(q_1q_2)^2 q_1 x_4(x_1x_2+x_3x_1+x_2x_3) =0,
%\end{align*}
%and similarly $\partial_4(\z{2}\z{1}+\z{1}\z{3}+\z{3}\z{2})=0$.
\epf
\begin{remark}
By  \cite{FK,MS} the following set is a basis of $\toba(\runow)$:
\begin{align*}
\{\z{1}^{n_{14}} \z{2}^{n_{24}} \z{3}^{n_{34}}\colon 0\leq n_{14}, n_{24}, n_{34}\leq 1 \} \cup  
\{\z{1}^{n_{14}} \z{2}\z{1} \z{3}^{n_{34}}\colon 0\leq n_{14}, n_{34}\leq 1 \}.
\end{align*}

\end{remark}

\subsection {Structure of $\rdosw=(\ad_cV_1)^2V_2$} \label{subsec:rdosw}
We set 
\begin{align}\label{eq:uij}
\hu{ij} &= \ad_cx_i (\z{j}) \in \rdosw,&  i, j&\in \Ib_3.
\end{align}

\begin{lemma}\label{lema:uij}
\begin{enumerate}[leftmargin=*,label=\rm{(\alph*)}]
\item\label{item:uij-props}
The subspace $\rdosw$ is a Yetter-Drinfeld submodule of $\toba(\nombre)$ with basis $\{\hu{12}, \hu{13}\}$. The action and the coaction satisfy
\begin{align}\label{eq:uij-trenza}
&\begin{aligned}
g_h \cdot \hu{1j} &=  q_1\omega^{(h+2)(1-j)} \hu{1 \, 5-j \,},\\ 
g_4 \cdot \hu{1j} &=- q_2^2\omega^2 \hu{1j}, \\
\delta (\hu{1j}) &= g_1g_jg_4 \ot \hu{1j},
\end{aligned}
&h \in \I_3,  j\in\I_{2,3}.
\end{align}
The braiding is of diagonal type, with matrix and Dynkin diagram
\begin{align*}
&\begin{pmatrix}
-1 & -\omega^2  \\
-\omega^2 & -1
\end{pmatrix},  & &\xymatrix{ \overset{-1}{\underset{\ }{\circ}}\ar  @{-}[rr]^{\omega}  &&
\overset{-1}{\underset{\ }{\circ}}}.
\end{align*}
Hence $\rdosw$ is of super type $\supera{\omega}{\I_2}$ \cite[Section 5.1]{AA}.

\item\label{item:uij-map} The inclusion induces an algebra map $\toba(\rdosw) \to \toba(\nombre)$. Hence the following relations hold:
\begin{align*} 
\hu{1j}^2&=0,\quad j\in \I_{2,3}, &
(\hu{12}\hu{13}+\omega^2\hu{13}\hu{12})^3&=0.
\end{align*}
\end{enumerate}
\end{lemma}

%The formula \eqref{eq:HV-link-xi-u1j} also holds for $j=1$ but we do not need it.

\pf \ref{item:uij-props}
Since the $x_i$'s are primitive in $\toba(\nombre)$, we have
\begin{align*}
\hu{ij} &=x_i\z{j}-m(g_i \cdot \z{j} \ot x_i) = x_i\z{j}+q_1\z{2i-j \,}x_i,& i,&j\in\I_3.
\end{align*}
As $x_i^2=0$, we have that $\hu{ii}=(\ad_c x_i)^2 x_4=0$. Hence $\rdosw$ is spanned by $\hu{ij}$, $i\ne j\in\I_3$.

\begin{claimo} If $i\neq j\in \I_3$, then 
\begin{align}\label{eq:uij-u1h}
\hu{ij} & = \omega^{(j-i)(1-i)}\hu{1 \, j-i+1 \, }
\end{align}
Hence $\rdosw$ is spanned by $\{\hu{12}, \hu{13}\}$.
\end{claimo}
We have to compute the skew derivations $\partial_h$ of both sides. If $h\neq i\in \I_3$, then $\partial_h(\hu{ij})=0$ by \eqref{eq:zi-skew}. For $h=i, 4$, 
\begin{align*}
\partial_i(\hu{ij}) &=g_i \cdot \z{j} + q_1 \z{2i-j \,}=0,
\\
\partial_4(\hu{ij})&=x_i \partial_4(\z{j})+q_1\partial_4(\z{2i-j \,})  g_4 \cdot x_i =-\omega^2x_ix_j+x_{2i-j}x_i.
\end{align*}
Hence \eqref{eq:uij-u1h} follows using \eqref{eq:rels-FK3}. Thus $\rdosw$ is spanned by $\{\hu{12}, \hu{13}\}$.
Now $\hu{12}$ and $\hu{13}$ are nonzero because $\partial_j \partial_4(\hu{1i})=-\omega^2x_i\neq0$, and linearly independent since $g_1g_ig_4 \neq g_1g_jg_4$ if $i\neq j$.
The action of $G$ and the braiding are obtained by direct computation.

\medbreak
Now \ref{item:uij-map} is a consequence of \cite[Theorem 2.6, Corollary 2.7(2)]{HS-JAlg} and \cite[Section 5.1]{AA}.
\epf

\begin{remark}
By \cite[Section 5.1.11]{AA}, $\toba(\rdosw)$ has a PBW basis:
\begin{align*}
 \{ \hu{13}^{n_{13}} (\hu{12}\hu{13}+\omega^2\hu{13}\hu{12})^{n_{1213}} \hu{12}^{n_{12}} \colon 0\leq n_{1213} < 3, 0\leq n_{12}, n_{13} <2  \}.
\end{align*}
\end{remark}

\subsection{Defining relations} \label{subsec:presentation-B(HV)}
Recall that 
\begin{align*}
\z{h} &= x_hx_4 - q_1  x_4 x_h, & 
\hu{ij} &= x_i\z{j}+q_1\z{2i-j \,}x_i, &
h, i, j & \in \I_3.
\end{align*}

The aim of this subsection is to give a minimal presentation of the Nichols algebra $\toba(\nombre)$. We start by stating the main result.

\begin{theorem} \label{thm:minimal-presentation}
The Nichols algebra $\toba(\nombre)$ is presented by generators $x_{i}$, $i \in \Ib_4$, and relations
\begin{align}
\label{eq:HV-xi^2}
x_{i}^{2}&=0,& &i  \in \Ib_3,
\\ \label{eq:HV-sum-xixj}
x_{1}x_{2} + x_{3}x_{1} +x_{2}x_{3} &=0,&  &x_{2}x_{1} + x_{1}x_{3} +x_{3}x_{2} =0,
\\ \label{eq:HV-x4^6}
x_4^6 &=0,
\\ \label{eq:HV-sum-u1iu1j}
(\hu{12}\hu{13}+\omega^2\hu{13}\hu{12})^3&=0,
\\ \label{eq:HV-ui1=u1k}
\hu{ i 1} - \omega \hu{1 \, 2 - i \,}  &= 0, & & i\in \I_{2,3},
\\ \label{eq:HV-link-x4-zh}
x_4\z{h}-q_2\z{h}x_4 &=0, &  & h  \in \I_3.
\end{align}
Let $\cG$ be the set of generators \eqref{eq:HV-xi^2}--\eqref{eq:HV-link-x4-zh}
of the ideal $\cJ(\nombre)$. Then $\cG$ is minimal.
The set
\begin{align}\label{D3-1:pbw}
x_4^{n_4}
\left\lbrace 
\begin{aligned}
\z{1}^{n_{14}}\z{2}^{n_{24}}\z{3}^{n_{34}}\\
\z{1}^{n_{14}}\z{2}\z{1}\z{3}^{n_{34}}
\end{aligned}
\right\rbrace 
\hu{12}^{n_{124}}\hu{12413}^{n_{124134}}\hu{13}^{n_{134}}
\left\lbrace 
\begin{aligned}
x_1^{n_1}x_2^{n_2}x_3^{n_3}\\
x_1^{n_1}x_2x_1x_3^{n_3}
\end{aligned}
\right\rbrace,
\end{align}
with $0\leq n_4 <6$, $0\leq n_{14}, n_{24}, n_{34}, n_{124}, n_{134}, n_1, n_2, n_3 <2$,  $0\leq n_{124134} <3$
is a basis of $\toba(\nombre)$, where $\hu{12413}=\hu{12}\hu{13}+\omega^2\hu{13}\hu{12}$.
\end{theorem}

Here the brackets mean that we choose either an element of the first line or else one of the second, with the restrictions for the $n_{\alpha}$.

\subsubsection{A first presentation of $\toba(\nombre)$}
We start with a set of defining relations of the Nichols algebra $\toba(\nombre)$. This set is far from being minimal but will help to prove Theorem \ref{thm:minimal-presentation}.

\medspace

By Remark \ref{rem:inclusion} and Lemmas \ref{lema:zi}, \ref{lema:uij}, 
the following relations hold in $\toba(\nombre)$:
\eqref{eq:HV-xi^2}, \eqref{eq:HV-sum-xixj}, \eqref{eq:HV-x4^6}, \eqref{eq:HV-sum-u1iu1j}, \eqref{eq:HV-ui1=u1k},
\begin{align} 
\label{eq:HV-zh^2}
\z{h}^{2}&=0,&  & h  \in \I_3,
\\ \label{eq:HV-sum-zizj}
\z{1}\z{2} + \z{3}\z{1} +\z{2}\z{3} &=0,&  &\z{2}\z{1} + \z{1}\z{3} +\z{3}\z{2} =0,
\\ \label{eq:HV-u1j^2}
\hu{1j}^2&=0,&  & j \in \I_{2,3},
\\ \label{eq:HV-uij=u1k}
\hu{ij}-\omega^{(j-i)(1-i)}\hu{1 \, j-i+1 \, }&=0 ,& &i\in \I_{2,3},  \, j\in\I_3, \, i\neq j,
\\ \label{eq:HV-link-xi-u1j}
x_i\hu{1j}-q_1\omega^{(i+2)(1-j)} \hu{1 \, 5-j \,}x_i&=0, & &i\in \I_3, \,  j\in\I_{2,3}.
\end{align}

\begin{lemma}\label{lemma:HV-presentation}
The algebra $\toba(\nombre)$ is presented by generators $x_{i}$, $i \in \Ib_4$, with relations \eqref{eq:HV-xi^2}--\eqref{eq:HV-link-xi-u1j}.
\end{lemma}

\pf  Let $\wtoba$ be the algebra presented by generators $\wx_{i}$, $i \in \Ib_4$, with relations \eqref{eq:HV-xi^2}--\eqref{eq:HV-link-xi-u1j}, 
replacing $\z{h}$ and $\hu{ij}$ by 
\begin{align*}
\wz{h}&=\wx_h \wx_4 - q_1\wx_4 \wx_h, & \wu{ij}&=\wx_i\wz{j}+q_1\wz{2i-j \,} \wx_i, && h, i, j \in \I_3.
\end{align*}

\begin{claim}\label{claim:exchange-rels}
There exists a surjective algebra map $\phi \colon \wtoba \to \toba(\nombre)$. 
\end{claim}
\pf
We already know that relations \eqref{eq:HV-xi^2}--\eqref{eq:HV-ui1=u1k},
and \eqref{eq:HV-zh^2}--\eqref{eq:HV-link-xi-u1j} hold in $\toba(\nombre)$. Hence it suffices to show that 
\eqref{eq:HV-link-x4-zh} also hold in $\toba(\nombre)$.
Notice that $\partial_i$, $i\in\I_3$, vanish both sides of \eqref{eq:HV-link-x4-zh}. Now
\begin{align*}
\partial_4(x_4\z{h})&=-q_2\omega^2\z{h}-\omega^2x_4x_h
&
\partial_4(\z{h}x_4)&=\omega^4x_hx_4+\z{h}.
\end{align*}
Hence
\begin{align*}
\partial_4(x_4\z{h} &-q_2\z{h}x_4)=-q_2\omega^2\z{h}-\omega^2x_4x_h- q_2
(\omega^4x_hx_4+\z{h})
\\
& =-q_2\omega^2 (x_hx_4-q_1x_4x_h) -\omega^2x_4x_h- q_2 (\omega x_hx_4+ x_hx_4-q_1x_4x_h)
\\
& =-q_2(1+\omega+\omega^2) x_hx_4 - (1+\omega+\omega^2)x_4x_h =0.
\end{align*}
so \eqref{eq:HV-link-x4-zh} holds in $\toba(\nombre)$.
\epf

\medbreak

Let $S_2 = \{\wx_4\}$, $S_{12}= \{\wz{1}, \wz{2}, \wz{3} \}$, $S_{112}= \{\wu{12}, \wu{13}\}$, $S_1= \{\wx_1, \wx_2, \wx_3\}$; 
let $\wtoba_{\alpha}$ be the subalgebra  of $\wtoba$ generated by $S_{\alpha}$ and
denote by $\phi_{\alpha}$ the restriction of $\phi$ to $\wtoba_{\alpha}$, $\alpha \in \varDelta_+$.

\begin{claim}\label{claim:dim-factores}
$\phi_{\alpha} \colon \wtoba_{\alpha} \to\toba(V_{\alpha})$ is an algebra isomorphism, $\alpha \in \varDelta_+$. 
\end{claim}
\pf
Since $S_{\alpha}$ generates $\wtoba_{\alpha}$ and is mapped by $\phi$ onto the set of generators of $\toba(V_{\alpha})$, it follows that $\phi (\wtoba_{\alpha}) = \toba(V_{\alpha})$.
Moreover, since \eqref{eq:HV-xi^2}--\eqref{eq:HV-link-xi-u1j} hold in  $\wtoba$, the defining relations of $\toba(V_{\alpha})$ are satisfied by the corresponding elements of $S_{\alpha}$. 
This yields an algebra map $\toba(V_{\alpha}) \to \wtoba_{\alpha}$ that is a section of $\phi_{\alpha}$.
\epf
\begin{claim}\label{claim:surjective}
The multiplication $\wtoba_{2} \ot \wtoba_{12} \ot \wtoba_{112} \ot \wtoba_1 \to \wtoba$ is surjective.
\end{claim}
\pf
We start with an auxiliar relation, that holds either in $\wtoba$ or $\toba(\nombre)$:
\begin{align}  \label{eq:HV-link-u1j-zh}
\hu{1j}\z{h}&=q_1\z{j+h-1 \,} \hu{1j}, & &h\in \I_3,  j\in\I_{2,3}.
\end{align}
In fact, by \eqref{eq:HV-uij=u1k} we may rewrite $\hu{1j}=\omega^m \hu{ih}$, where $i=h-j+1$ and $m\in\I_{0,2}$ depends on $h$ and $j$. Using \eqref{eq:HV-zh^2} twice,
\begin{align*}
\hu{1j} \z{h}&= \omega^m \hu{ih} \z{h} = \omega^m 
(x_{i} \z{h} + q_1 \z{j+h-1 \,} x_{i} ) \z{h}
\\
& = \omega^m q_1 \z{j+h-1 \,} ( \hu{i h} - q_1 \z{j+h-1 \,} x_{i} ) = \omega^m q_1 \z{j+h-1 \,}  \hu{ih}
\\
& = q_1 \z{j+h-1 \,}  \hu{1j}.
\end{align*}

Let $\mD_2=\wtoba_{2}$, $\mD_{12}=\wtoba_{2} \wtoba_{12}$,
$\mD_{112}=\wtoba_{2} \wtoba_{12} \wtoba_{112}$,
$\mD_1=\wtoba_{2} \wtoba_{12} \wtoba_{112} \wtoba_1$. 
Our goal is to prove that $\wtoba=\mD_1$.
Since $\mD_1$ contains $1$, this can be achieved showing that $\mD_1$ is a left ideal of $\wtoba$, which reduces to prove that $\wx_i \mD_1 \subset \mD_1$ for all $i \in \I_4$: The case $i=4$ is straightforward.
We proceed by steps. 

\begin{step}\label{step:surjective1}
 $\mD_{112} S_{12} \subset \mD_{112}$, $\mD_{1} S_{112} \subset \mD_{1}$. 
\end{step}

We start with the first inclusion. 
We have to prove that $\mD_{112} \wz{h} \subset \mD_{112}$ for all $h\in\I_{3}$. The subspace $\mD_{112}$ is spanned by monomials 
\begin{align*}
y&=\wx_4 ^{n} \wz{h_1} \cdots \wz{h_k} \wu{1 j_1}\cdots\wu{1j_l}, & n, \, & k, l \geq 0, \, h_i\in\I_{3}, \, j_t \in\I_{2,3}.
\end{align*}
Fix $h\in\I_{3}$. If $l=0$, then $y \wz{h} \in \mD_{112}$. If $l>0$, then
\begin{align*}
y \wz{h} &= q_1^l \wx_4^{n} \wz{h_1} \cdots \wz{h_k} \wz{h'}
\wu{1j_1}\cdots\wu{1j_l} & \text{for some }& h'\in\I_{3},
\end{align*}
by repeated applications of \eqref{eq:HV-link-u1j-zh}. Thus $y \wz{h} \in \mD_{112}$.

Now we notice that $\mD_{1}$ is spanned by elements
\begin{align*}
z&=\wx_4^{n} \wz{h_1} \cdots \wz{h_k} \wu{1j_1}\cdots\wu{1j_l} \wx_{i_1} \cdots \wx_{i_m},
& n, & k, l \geq 0, \, h_r, i_s\in\I_{3}, \, j_t \in\I_{2,3}.
\end{align*}
If $m=0$, then $y \wu{1j} \in \mD_{1}$ by definition. If $m>0$, then
\begin{align*}
z \wu{1j} &= q_1^m \omega^b \wx_4^{n} \wz{h_1} \cdots \wz{h_k} \wu{1,j_1}\cdots\wu{1,j_l}\wu{1j'} \wx_{i_1} \cdots \wx_{i_m}
\end{align*}
for some $b \in\N_0$, $j' \in \I_{2,3}$, 
by \eqref{eq:HV-link-xi-u1j}. Hence $y \wu{1j} \in \mD_{1}$.

\begin{step}\label{step:surjective2}
For all $i \in \I_3$,
\begin{align*}
\ku \langle \wx_i \rangle \wtoba_2 &\subset \mD_{12} \, \ku \langle \wx_i \rangle, &
\ku \langle \wx_i \rangle \wtoba_{12} &\subset \mD_{112} \, \ku \langle \wx_i \rangle, &
\ku \langle \wx_i \rangle \wtoba_{112} &\subset \mD_{1} \, \ku \langle \wx_i \rangle.
\end{align*}
\end{step}

For the first inclusion, we only need to show that $\wx_i \wtoba_2 \subset \mD_{12} \, \ku \langle \wx_i \rangle $ since $\wx_i^2=0$ and $\wtoba_2 \subset \mD_{12}$.
The claim holds by the following formula, which follows inductively from \eqref{eq:HV-link-x4-zh}:
\begin{align*}
\wx_i \wx_4^k&=\sum_{h=0}^{k-1} (q_2^{h-(k-1) }q_1^h) \wx_4 ^{k-1} \wz{i}+ q_1^ k \wx_4^ k \wx_i, & i \in \I_3, k\geq 1.
\end{align*}

For the second, it suffices to prove that $\wx_i \wtoba_{12} \subset \mD_{112} \, \ku \langle \wx_i \rangle $.  We prove by induction on $k$ that $\wx_i \wz{h_1} \cdots \wz{h_k} \in \mD_{112} \, \ku\langle \wx_i \rangle$. 
The case $k=0$ is trivial, while $k=1$ is just the definition of the $\wu{ij}$'s. We assume it holds for $k$. Then, for any $h \in \I_3$, 
\begin{align*}
\wx_i \wz{h_1} \cdots \wz{h_k} \wz{h} & \in  \mD_{112} \, \ku\langle \wx_i \rangle \wz{h} 
\subset  \mD_{112} \, (S_{112} + S_{12} \ku \langle \wx_i \rangle) \subset \mD_{112} \ku \langle \wx_i \rangle,
\end{align*}
where we use inductive hypothesis, Step \ref{step:surjective1} and the following facts:
\begin{align*}
\ku \langle \wx_i \rangle \wz{h} & \subset  S_{112} +  S_{12} \ku \langle \wx_i \rangle, & \mD_{112} S_{112} & \subset \mD_{112}.
\end{align*}

The last inclusion follows from \eqref{eq:HV-link-xi-u1j}.

\begin{step}\label{step:surjective3}
$\mD_2 \mD_{12} \subset \mD_{12}$, $\mD_{12} \mD_{112} \subset \mD_{112}$.
\end{step}
 
The first inclusion can be written as $\wtoba_2 \wtoba_2 \wtoba_{12} \subset \wtoba_2 \wtoba_{12}$, and the proof follows since $\wtoba_2$ is a subalgebra. The proof of the second inclusion is similar since $\wtoba_{12}\wtoba_2  =\wtoba_2 \wtoba_{12}$, and $\wtoba_{12}$, $\wtoba_2$ are subalgebras.

\begin{step}\label{step:surjective4}
$\wx_i \mD_{112} \subset \mD_{112}\, \ku \langle \wx_i \rangle$, $i \in \I_3$.
\end{step}
 
By \eqref{eq:HV-link-xi-u1j}, $\ku \langle \wx_i \rangle \wtoba_{112} = \wtoba_{112}  \ku \langle \wx_i \rangle$.
Now we have
\begin{align*}
\wx_i \mD_{112} &= \wx_i \wtoba_2 \wtoba_{12} \wtoba_{112}
\subset \mD_{12} \ku \langle \wx_i \rangle \wtoba_{12} \wtoba_{112} 
\subset \mD_{12} \mD_{112} \ku \langle \wx_i \rangle \wtoba_{112}
\\
& \subset \mD_{112} \wtoba_{112} \ku \langle \wx_i \rangle = \mD_{112} \ku \langle \wx_i \rangle,
\end{align*}
where the first two inclusions follow from Step \ref{step:surjective2}, the third from Step \ref{step:surjective3} and
$ \ku \langle \wx_i \rangle \wtoba_2 \subset \wtoba_2  \ku \langle \wx_i \rangle$,
and the last equality from $\mD_{112} \wtoba_{112} = \mD_{112}$.

\begin{step}\label{step:surjective5}
$\wx_i \mD_{1} \subset \mD_{1}$, $i \in \I_3$.
\end{step}

Indeed, by Step \ref{step:surjective4} and the fact $\ku \langle \wx_i \rangle \wtoba_1 = \wtoba_1$ it follows
\begin{align*}
\wx_i \mD_{1} = \wx_i \mD_{112} \wtoba_1 
\subset \mD_{112} \ku \langle \wx_i \rangle \wtoba_1 
= \mD_{112} \wtoba_1 = \mD_{1}.
\end{align*}
This concludes the proof of Claim \ref{claim:surjective}. 
\epf

By Claim \ref{claim:surjective} in $\ast$,  Claim \ref{claim:dim-factores} and Theorem \ref{thm:D3-1} in $\star$ and
Claim \ref{claim:exchange-rels} in $\dagger$, we have 
\begin{align*}
\dim\wtoba \overset{\ast}{\leq} \dim\wtoba_2\dim\wtoba_{12}\dim\wtoba_{112}\dim\wtoba_1 \overset{\star}{=} \dim\toba(\nombre) \overset{\dagger}{\leq} \dim\wtoba.
\end{align*}
Since $\phi$ is surjective it must be an isomorphism.
\epf

\subsubsection{Proof of Theorem \ref{thm:minimal-presentation}}
The statement about the basis follows from Theorem \ref{thm:D3-1} and Lemmas \ref{lema:zi}, \ref{lema:uij}.

\medbreak

We seek now for a minimal set of generators of the ideal $\cJ ( \nombre )$. 
Hence we need to obtain relations in Lemma \ref{lemma:HV-presentation} from relations in Theorem \ref{thm:minimal-presentation}.

\begin{lemma}\label{lemma:discard-rels}
Let $S$ be a quotient algebra of $T ( \nombre )$
such that \eqref{eq:HV-xi^2}, 
\eqref{eq:HV-sum-xixj}, \eqref{eq:HV-ui1=u1k} and \eqref{eq:HV-link-x4-zh} 
hold in $S$.
Then \eqref{eq:HV-zh^2}, \eqref{eq:HV-sum-zizj}, \eqref{eq:HV-uij=u1k}, \eqref{eq:HV-link-xi-u1j} and \eqref{eq:HV-u1j^2} also hold.

%Let $S$ be a quotient algebra of $T ( \nombre )$.
%\begin{enumerate}[leftmargin=*,label=\rm{(\alph*)}]
%\item\label{item:discard-zh^2} If \eqref{eq:HV-xi^2} and \eqref{eq:HV-link-x4-zh} hold in $S$, then \eqref{eq:HV-zh^2} also holds.
% \item\label{item:discard-sum-zizj} If \eqref{eq:HV-sum-xixj} and \eqref{eq:HV-link-x4-zh} hold in $S$, then \eqref{eq:HV-sum-zizj} also holds.
% \item\label{item:discard-u23} If \eqref{eq:HV-sum-xixj} and \eqref{eq:HV-ui1=u1k} hold in $S$, then \eqref{eq:HV-uij=u1k} also holds.
% \item\label{item:discard-link-xi-u1j} If \eqref{eq:HV-xi^2} and \eqref{eq:HV-uij=u1k} hold in $S$, then \eqref{eq:HV-link-xi-u1j} also holds.
% \item\label{item:discard-link-u1j-zh} 
%If \eqref{eq:HV-zh^2} and \eqref{eq:HV-uij=u1k} hold in $S$, then \eqref{eq:HV-link-u1j-zh} also holds.
% \item\label{item:discard-u1j^2} If \eqref{eq:HV-xi^2} and \eqref{eq:HV-link-xi-u1j} hold in $S$, then \eqref{eq:HV-u1j^2} also holds.
%\end{enumerate}
\end{lemma}

\pf
Let $h, j \in \I_3$. As we assume that \eqref{eq:HV-link-x4-zh} holds, untwining  \eqref{eq:HV-link-x4-zh} we get
$x_4 (x_hx_4 - q_1 x_4 x_h) = q_2 (x_h x_4 -q_1 x_4 x_h) x_4$, which can be rewritten as 
\begin{align*}
(1-\omega) x_4x_h x_4 &= q_2 x_h x_4^2 + q_1 x_4^2 x_h.
\end{align*}
Using this equality we compute
\begin{align}\label{eq:zhzj}
\begin{aligned}
(1-\omega) \z{h} \z{j} &= (1-\omega)(x_hx_4 - q_1 x_4 x_h)(x_jx_4 - q_1 x_4 x_j)\\
&=x_h (q_2 x_j x_4^2 + q_1  x_4^2 x_j ) - q_1 (1 - \omega) x_h x_4^2 x_j - \\
& \qquad- q_1 (1 - \omega) (x_4 x_h x_j x_4) + q_1^2  (q_2 x_h x_4^2 + q_1  x_4^2 x_h ) x_j \\
%=q_2 x_h x_j x_4^2 + q_1 (1 - (1 - \omega) + q_1 q_2) x_h x_4^2 x_j - q_1 (1 - \omega) (x_4 x_h x_j x_4) + q_1^3 x_4^2 x_h x_j \\ 
&= q_2 x_h x_j x_4^2  - q_1 (1 - \omega) (x_4 x_h x_j x_4) + q_1^3 x_4^2 x_h x_j.
\end{aligned}
\end{align}
Now \eqref{eq:HV-zh^2} follows from \eqref{eq:zhzj} for $j=h$ and \eqref{eq:HV-xi^2}.
%\begin{align*}
%(1-\omega) \z{h}^2&= q_2 x_h^2 x_4^2  - q_1 (1 - \omega) (x_4 x_h^2 x_4) + q_1^3 x_4^2 x_h^2 = 0, & h \in \I_3.
%\end{align*}

For \eqref{eq:HV-sum-zizj}, let $i \in \I_{2,3}$. Using \eqref{eq:zhzj} again and \eqref{eq:HV-sum-xixj},
\begin{align*}
\z{1} \z{i} + & \z{5-i \,} \z{1} + \z{i} \z{5-i \,} = q_2 (x_1 x_i + x_{5-i} x_1 + x_i x_{5-i} ) x_4 ^2 - \\
& - q_1 (1 - \omega) x_4 ( x_1 x_i + x_{5-i} x_1 + x_i x_{5-i} )x_4 
\\
& + q_1^3 x_4^2 (x_1 x_i + x_{5-i} x_1 + x_i x_{5-i} ) = 0.
\end{align*}
Next we prove \eqref{eq:HV-uij=u1k}. Let $i \in \I_{2, 3}$. By \eqref{eq:HV-sum-xixj}, $\hu{1 i} + \hu{5 - i\, 1} + \hu{i \, 5 - i \,} = 0$, so
\begin{align*}
\hu{i \, 5 - i \,} = - \hu{1 i} - \hu{5-i \, 1} \overset{\eqref{eq:HV-ui1=u1k}}{=} - (1+\omega) \hu{ 1 i } = \omega^2 \hu{ 1 i}.
\end{align*} 

Now we turn to \eqref{eq:HV-link-xi-u1j}. Let $i, j \in \I_3$. Since $x_i^2 = 0$ we have that
\begin{align} \label{eq:HV-link-xi-uij}
x_i \hu{ij} 
%= x_i(x_i\z{j}+q_1 \z{2i-j}x_i) = q_1(x_i \z{2i-j}+q_1 \z{j} x_i)x_i 
= q_1 \hu{i \, 2i - j \,} x_i .
\end{align}
We now split the proof of the six relations of \eqref{eq:HV-link-xi-u1j} in different cases. Assume first $i = 1$. Then \eqref{eq:HV-link-xi-u1j} means 
$x_1 \hu{1j} =  q_1 \hu{1 \, 2 - j\,} x_1$ for any $j \in \I_3$, which is \eqref{eq:HV-link-xi-uij}. 
Now assume $i \in \I_{2, 3}$. By \eqref{eq:HV-uij=u1k} and \eqref{eq:HV-link-xi-uij} we get
\begin{align*}
x_i \hu{1 i} &= \omega  x_i \hu{i \, 5- i \,} =q_1  \omega   \hu{i1} x_i = q_1 \omega ^2 \hu{ 1 \, 5- i \,} x_i ,\\ 
x_i \hu{1 \, 5-i \,} &= \omega^2 x_ i \hu{i1} =q_1 \omega^2 \hu{i \, 5-i \,} x_i = q_1 \omega \hu{1i} x_i.
\end{align*}

% \ref{item:discard-link-u1j-zh} Let $h \in \I_3$ and $j \in \I_{2,3}$. By \eqref{eq:HV-uij=u1k}, we may rewrite $\hu{1\, j}=\omega^m \hu{i\, h}$ with $i=h-j+1$ and some integer $m$. Now compute
%\begin{align*}
%\hu{1\, j} \z{h}&= \omega^m \hu{h-j+1,h} \z{h} = \omega^m (x_{h-j+1} \z{h} + q_1 \z{j+h-1 \,} x_{h-j+1} ) \z{h} =  \\
% & = \omega^m q_1 \z{j+h-1 \,} ( \hu{h-j+1, h} - q_1 \z{j+h-1 \,} x_{h-j+1} ) = \omega^m q_1 \z{j+h-1 \,}  \hu{h-j+1, h}\\
% & = q_1 \z{j+h-1 \,}  \hu{1, j}.
%\end{align*}

Finally we consider \eqref{eq:HV-u1j^2}. 
By \eqref{eq:HV-xi^2}, $\hu{11}= (\ad_c x_1^2)x_4=0$. Hence
\begin{align*}
\hu{1j}^2 &= (x_1 \z{j} + q_1 \z{2j-1 \,} x_1) \hu{1j} = q_1^{-1}  x_1 \hu{1j} \z{1} + q_1^2 \z{5-j \,} \hu{1 \, 5-j \,} x_1 \\
&= \hu{1 \, 5-j \,} x_1 \z{1} + q_1 \hu{1 \, 5-j \,} \z{1} x_1 = \hu{1 \, 5-j \,} \hu{11}= 0,
\end{align*}
for all $j \in \I_{2,3}$.
\epf

We proceed with the proof of the statement about the presentation of $\cJ(\nombre)$. By Lemmas \ref{lemma:HV-presentation} and \ref{lemma:discard-rels} $\cG$ generates $\cJ (\nombre)$ as an ideal.
Hence it remains to show the minimality: It is enough to show that if $r \in \cG$, then the ideal $\cI_r$ generated by $\cG \setminus r$ is not the whole $\cJ ( \nombre )$.
The tensor algebra $T ( \nombre )$ has a $\N_0^2$-grading with 
$x_1, x_2, x_3$ sitting on degree $(1,0)$ and $x_4$ on degree $(0,1)$. 
The relations \eqref{eq:HV-xi^2}--\eqref{eq:HV-link-x4-zh} are homogeneous and the $\Bbbk\Gamma$-coaction is of the shape $\delta(r)=g_r\ot r$ for some $g_r\in\Gamma$.
The next table contains $\N_0^2$-degrees and $g_r$ for all $r\in\cG$.

\medbreak

\begin{tabular}{ l r l c c}
\, & Relation & \, & $\N_0^2$-degree & $g_r$ \\
\hline
\eqref{eq:HV-xi^2}   & $x_{i}^{2}$,    &  $i  \in \Ib_3$,       & $(2,0)$  & $g_i^2$
\\
\eqref{eq:HV-sum-xixj}   & $x_1 x_i + x_{5-i} x_1 + x_i x_{5-i} $,   & $i  \in \I_{2, 3}$, & $(2,0)$  & $g_1 g_i$ 
\\
\eqref{eq:HV-x4^6} & $x_4^6$,      & \,    & $(0,6)$  & $g_4^6$
\\
\eqref{eq:HV-sum-u1iu1j}  & $(\hu{12}\hu{13}+\omega^2\hu{13}\hu{12})^3$,  & \,   & $(12,6)$ & $g_1 ^ {12} g_4 ^ 6$ 
\\
\eqref{eq:HV-ui1=u1k} & $\hu{i1}-\omega \hu{1\,5-i\,}$,    & $i\in \I_{2,3}$,    &$(2,1)$    & $g_1 g_{5-i} g_4$ 
\\
\eqref{eq:HV-link-x4-zh}   & $x_4\z{h} -q_2\z{h}x_4$,   &  $h  \in \I_3$ ,  & $(1,2)$   & $g_h g_4^2$
\end{tabular}

\medbreak

Consider on $\N_0^2$ the partial order $(a,b) \preceq (a', b')$ if and only if $ a\leq a'$ and $b\leq b'$. 
This induces an order on the set of homogeneous elements of $T (\nombre)$.

If $r$ is \eqref{eq:HV-x4^6} then it is minimal in $\cG$, so $r \notin \cI_r$ since its $\N_0^2$-degree is minimal.
Now fix $r$ one of the words with degree $(2,0)$. Since no element in $\cG$ has degree less than $(2,0)$, the only way to have 
$r \in \cI_r$ is to write it as a linear combination of the other elements in  $\cG$ with degree $(2,0)$. 
This is impossible since the $g_r$ are pairwise different so they are linearly independent. The same argument holds if $r$ is one of  \eqref{eq:HV-link-x4-zh}.

Now we fix a relation $r_i = \hu{i1}-\omega \hu{1 \, 5-i \,}$ in \eqref{eq:HV-ui1=u1k}. Suppose that $r_i \in \cI_r$: By $\N_0^2$-degree restrictions we must have
\begin{align*}
r_i &= \nu \, r_{5-i} + \sum_{r\in\cG, \deg r= (2,0)} \lambda_r \,  x_4r + \mu_r \,  r x_4, & \nu, \lambda_r, \mu_r & \in \Bbbk.
\end{align*}
Now we look at the $\Bbbk\Gamma$-coaction. Since $g_{r_i}=g_1 g_{5-i} g_4$, it follows that $\nu = 0$ and $\mu_r =\lambda_r=0$ if $g_r\neq g_1 g_{5-i}$. So our sum becomes
\begin{align} \label{eq:HV-ui1=u1k-minimal}
r_i= \lambda x_4 (x_1  x_{5-i} + x _ i x_1 +  x_{5-i} x_i) + \mu  (x_1  x_{5-i} + x _ i x_1 +  x_{5-i} x_i) x_4.
\end{align}
But untwining the definition we get 
\begin{align*}
r_i &= \hu{i1}-\omega \hu{1\, 5-i \,} = ( - \omega x_1 x_{5-i} + x_i x_1) x_4 + q_1^2 x_4 (\omega x_i x_1 - x_{5-i} x_i)
\\ 
& \qquad  + q_1 ( \omega x_1 x_4 x_ {5-i} + \omega ^ 2 x_i x_4 x_1 + x_{5-i} x_4 x_i),
\end{align*}
and the term of the form  $x_1 x_4 x_{5-i}$ does not appear on the right-hand side of \eqref{eq:HV-ui1=u1k-minimal}. This contradiction shows that $r\notin \cI_r$.

Finally, we show using GAP that \eqref{eq:HV-sum-u1iu1j} is not in the ideal generated by the other relations. 
\qed

\bigbreak

\section{Pre-Nichols algebras of $\nombre$}\label{sec:preNichols}

We have two purposes in this \S. On the one hand, we study finite-dimensional pre-Nichols algebras of $\nombre$ to conclude that any pointed Hopf algebra with infinitesimal braiding $\nombre$ is generated by skew-primitive and group-like elements. On the other hand, we introduce a pre-Nichols algebra $\widetilde{\toba}(\nombre)$ of $\nombre$ which plays the role of distinguished pre-Nichols algebras for braidings of diagonal type \cite{An-distinguished}: 
the Gelfand Kirillov dimension of $\widetilde{\toba}(\nombre)$ is finite and $\widetilde{\toba}(\nombre)$ has a skew central Hopf subalgebra $\cZ(\nombre)$ such that $\toba(\nombre)$ is the quotient of $\widetilde{\toba}(\nombre)$ by the ideal generated by $\cZ(\nombre)$.

\subsection{Generation in degree one}
Here we study finite-dimensional pre-Nichols algebras $\wtoba$ of $\nombre$. 
We recall that there exists a homomorphism $T(\nombre) \twoheadrightarrow \wtoba$ of graded Hopf algebras in $\ydg$.
In the following lemmas we proceed as in \cite[Lemma 5.4]{AS-Annals} to show that this map factors through some relations of the defining ideal $\J (\nombre)$. We start by looking at those relations which are primitive in $T(\nombre)$.

\begin{rem}\label{rem:comultiplication-relations}
By direct computation, 
\begin{align*}
x_i^2, \, x_1 x_i + x_{5-i} x_1 + x_i x_{5-i}, \, x_4^6 &\in\mathcal{P}(T(\nombre)).
\end{align*}
For \eqref{eq:HV-link-x4-zh} and \eqref{eq:HV-ui1=u1k}, the comultiplication in $T(\nombre)$ is
\begin{align*} 
\Delta (x_4 \z{h} - q_2 \z{h} x_4 ) &= (x_4 \z{h} - q_2 \z{h} x_4 ) \ot 1 + 1 \ot (x_4 \z{h} - q_2 \z{h} x_4 )
, \\
\Delta (\hu{i1} - \omega \hu{1 \,5-i\,}) &= (\hu{i1} - \omega \hu{1\,5-i\,}) \ot 1 + 1 \ot (\hu{i1} - \omega \hu{1\,5-i\,})  \\
&\qquad + (x_1 x_{5-i} + x_i x_1 + x_{5-i} x_i) \ot x_4.
\end{align*}
Indeed, these formulas follow from the following ones:
\begin{align*} 
\Delta ( \z{h} ) &= \z{h} \ot 1 + 1 \ot \z{h} - \omega^2 x_h \ot x_4,
\\
\Delta ( \hu{ij} ) &= \hu{ij} \ot 1 + 1 \ot \hu{ij}  + x_i \ot \z{j} + \omega x_j \ot \z{2i-j \,} + \omega^2 x_{2i-j} \ot \z{i}  \\
&\qquad + ( x_{2i-j} x_i - \omega ^2 x_i x_j ) \ot x_4. 
\end{align*}
\end{rem}

\begin{lemma}\label{lemma:S-discard-G0-G1}
Let $\wtoba$ be a finite-dimensional pre-Nichols algebra of $\nombre$. Then \eqref{eq:HV-xi^2}, \eqref{eq:HV-sum-xixj}, \eqref{eq:HV-x4^6} and \eqref{eq:HV-link-x4-zh} hold in $\wtoba$.
\end{lemma}
\pf
First we note that the Nichols algebra of the primitive elements $\mP (\wtoba)$ is finite dimensional. 
Indeed, $\mP (\wtoba) \# \ku \Gamma$ is contained in the first term of the coradical filtration of $\wtoba \# \ku \Gamma$, see \cite[Lemma 5.4]{AS-Annals}.

Since $T(\nombre) \to \wtoba$  is a homomorphism of braided Hopf algebras, the elements
\eqref{eq:HV-xi^2}, \eqref{eq:HV-sum-xixj}, \eqref{eq:HV-x4^6} and \eqref{eq:HV-link-x4-zh}
are primitive in $\wtoba$.  The strategy now is to build braided subspaces
of $\mP(\wtoba)$ which are either zero or generate an infinite-dimensional Nichols algebra.
By direct computation, if $r$ is one of the relations in 
\eqref{eq:HV-xi^2}, \eqref{eq:HV-sum-xixj} or \eqref{eq:HV-x4^6}, then $c( r \ot r) = r \ot r$. Hence $r=0$.

Now we turn to \eqref{eq:HV-link-x4-zh}. Fix $h \in \I_3$ and set $r_h = x_4 \z{h} - q_2 \z{h} x_4  \in \mP (\wtoba)$.
Since $g_4$ is central, $r_h$ is homogeneous of degree $g_h g_ 4 ^ 2$. 
Suppose that $r_h$ is non-zero, we may consider the $2$-dimensional subspace $W_h = \ku x _h \oplus \ku r_h \subset \mP (\wtoba)$. 
The braiding of $W_h$ is diagonal, with braiding matrix and Dynkin diagram
\begin{align*}
&\begin{pmatrix}
-1 & -q_1^2  \\
-q_2^2 & - \omega
\end{pmatrix},  & &\xymatrix{ \overset{-1}{\underset{\ }{\circ}}\ar  @{-}[rr]^{\omega ^ 2}  &&
\overset{- \omega}{\underset{\ }{\circ}}},
& & - \omega \in \Gp_6.
\end{align*}
%We compute the braiding: 
%\begin{align*}
%c( x_h \ot y_h) &= g_h \cdot y_h \ot x_h  = - q_1 ^ 2  y_ h \ot x_h, \\
%c( y_h \ot x_h) &= g_h g_4 ^ 2 \cdot x_h \ot y_h  = - q_ 2 ^ 2  x_h \ot y_h, \\
%c( y_h \ot y_h) &= g_h g_4 ^ 2 \cdot y_h \ot y_h  = - \omega  y_ h \ot y_ h.
%\end{align*}
This diagram does not appear in \cite[Table 1]{H-rank2root}. Hence $\dim \toba(W_h)=\infty$  and we have a contradiction since $\toba(W_h)\hookrightarrow \toba( \mP (\wtoba))$. Thus $r_h = 0$.
\epf

\begin{lemma} \label{lemma:S-discard-G2}
Let $\wtoba$ be as above. Then \eqref{eq:HV-ui1=u1k} holds in $\wtoba$.
\end{lemma}
\pf
Set $y_1 = \hu{31} - \omega \hu{12}$, $y_2 = \hu{21} - \omega \hu{13}$ and $W = \ku y_1 + \ku y_2$; 
we have $W \subset  \mP(\wtoba)$ by Lemma \ref{lemma:S-discard-G0-G1}. 
The $\Bbbk\Gamma$-coaction is given by $g_1g_2g_4$ and $g_1g_3g_4$ respectively. Let us compute the $\Gamma$-module structure of $W$.
As \eqref{eq:HV-sum-xixj} holds in $\wtoba$ by the previous step, $ \hu{1i} + \hu{2- i \, 1} + \hu{i \, 2- i \,}= 0$, $i \in \I_{2, 3}$. Then
\begin{align}\label{eq:action-ui1=u1k}\begin{aligned}
g_1 \cdot y_1 &=  q_1 (\hu{21} - \omega \hu{13} ) = q_1 y_2,\\
g_ 2 \cdot y_1 &= q_1 ( \hu{13} - \omega \hu{32} )=   q_1 ( (1 + \omega ) \hu{13} + \omega \hu{21} ) = q_1 \omega y _ 2, \\
g_ 3 \cdot y_1 &= q_1 ( \hu{32} - \omega \hu{21} ) = q_1 ( - \hu{13} + (- 1 -\omega )  \hu{21} ) = q_1 \omega ^ 2  y _ 2, \\
g_4 \cdot y_1 & = - q_2 ^ 2 \omega ^ 2 y _1,
\end{aligned}\end{align}
and similarly 
\begin{align*}
g_1 \cdot y_2  = q_1 y _1, \quad g_ 2 \cdot y_2 =q_1 \omega^2 y _ 1, \quad g_ 3 \cdot y _ 2 = q_1 \omega  y _ 1, 
\quad g_4 \cdot y _ 2  = - q _ 2 ^ 2 \omega ^ 2 y _2.
\end{align*}
Assume $W \neq 0 $. Since $g_1$ permutes the generators, we have  $y_1, y _2 \neq 0$; moreover they are linearly independent since they have different $\Gamma$-degrees.
A straightforward computation shows that $W$ has not $\Gamma$-stable $1$-dimensional subspace, so $W$ is a simple Yetter-Drinfeld module over $\Gamma$. 

As $\supp V\oplus W$ generates the subgroup $\Gamma'$ generated by $g_i$, $i\in\I_4$, we may realize our Nichols and pre-Nichols algebras over $\Gamma'$. Hence we may (and will) assume that $\Gamma=\Gamma'$.
We may now evoke the classification theorem on \cite{HV-rank2}. Indeed, consider the pair $(V, W)$ of simple objects in $\ydg$. 
As $W\subset \wtoba^2$, we have $V \cap W = 0$. By Remark \ref{rem:comultiplication-relations} and Lemma \ref{lemma:S-discard-G0-G1}, we have $ V \oplus W \subset \mP (\wtoba)$, so
$\dim \toba(V \oplus W)<\infty$. Thus
$V \oplus W$ is one of the five-dimensional braided vector spaces in \cite[Theorem 2.1]{HV-rank2}; moreover it is one of \cite[Examples 1.9--1.11]{HV-rank2} since $V$ is the quandle of transpositions of $\mathbb{S}_3$. Next we compute the action of $\Gamma$ on $V \oplus W$:

\medbreak

\begin{tabular}{ L | C C C C C }
\, & x_1 & x_2 & x_3 & y_1 & y_2 \\
\hline
g_1 & -x_1 & -x_3 & -x_2 & q_1 y_2 & q_1 y_1\\
g_2 & -x_3 & -x_2 & -x_1 & q_1 \omega y_2 & q_1\omega ^ 2  y_1\\
g_3 & -x_2 & -x_1 & -x_3 & q_1 \omega ^ 2 y_2 & q_1\omega  y_1\\
g_1g_2g_4 & q_2 x_2 & q_2x_3 & q_2 x_1 & - \omega^2 y_1 & - y_2\\
g_1g_3g_4 & q_2 x_3 & q_2x_1 & q_2 x_2 & -  y_1 & - \omega^2 y_2
\end{tabular}

\bigbreak

Hence $c_{W,V} c_{V, W} \neq \id_{V\ot W}$
and the braiding of $W$ is diagonal, with matrix $\begin{pmatrix}
-\omega ^ 2 & -1  \\
-1 & - \omega^2
\end{pmatrix}$, $-\omega^2 \in \Gp_6$. But all the $W$'s in \cite[Examples 1.9 - 1.11] {HV-rank2} have vertices labeled with $-1$, which is a contradiction. Thus $W=0$.
\epf

%Since we did not explicitly compute its comultiplication on $T(\nombre)$, we will need the following result.
%
%\begin{lemma*} \cite[Lemma 3.2]{An-diagonal} 
%Let $V \in \ydh$ and $ I \subset T(V)$ a braided homogeneous biideal of $T(V)$, so there is a surjective morphism of braided graded Hopf algebras 
%$\pi: R := T(V) / I \twoheadrightarrow \toba (V)$. Let $x \in \ker \pi$, $x \neq 0$ of minimal degree $ \geq 2$. Then $x$ is primitive.
%\end{lemma*} 

We can now prove the main result of this section, which states that $\toba (\nombre)$ is the unique finite dimensional pre-Nichols, respectively post-Nichols, algebra of $\nombre$.
In particular, we will show that the top degree relation \eqref{eq:HV-sum-u1iu1j} holds in any finite dimensional pre-Nichols algebra.

\begin{theorem}\label{thm:pre-post-Nichols}
\begin{enumerate}[leftmargin=*,label=\rm{(\alph*)}]
\item\label{item:prenichols} Let $\wtoba =\oplus_{n\geq0} \wtoba^n$ be a finite-dimensional pre-Nichols algebra of $\nombre$. Then $\wtoba\simeq\toba(\nombre)$.

\item\label{item:postnichols} Let $\cL =\oplus_{n\geq0} \cL^n$ be a finite-dimensional post-Nichols algebra of $\nombre$. Then $\cL\simeq\toba(\nombre)$.
\end{enumerate}
\end{theorem}

\pf
For \ref{item:prenichols} we proceed as in \cite [Theorem 4.1]{An-diagonal}. By definition, we may identify $\wtoba= T(\nombre)/ I$, where $I$ is a 
graded Yetter-Drinfeld submodule and Hopf ideal of $T(\nombre)$ generated by homogeneous elements of degree $\geq 2$, and $I \subseteq \cJ (\nombre)$. Let
 $\pi \colon \wtoba \twoheadrightarrow \toba (\nombre)$ be the canonical projection of graded braided Hopf algebras.
Assume $\cJ (\nombre) \supsetneq I $, hence one of the generators in Theorem \ref{thm:minimal-presentation} does not belong to $I$. 
By Lemmas \ref{lemma:S-discard-G0-G1} and \ref{lemma:S-discard-G2}, it must be $r = (\hu{12}\hu{13}+\omega^2\hu{13}\hu{12})^3 \notin I$.
As $\cG$ is a minimal set of defining relations of $\cJ(\nombre)$ and this ideal is graded, $r$ has minimal degree $\geq 2$ among 
non-trivial elements in $\ker\pi$.
Then $r$ is primitive in $\wtoba$ by \cite [Lemma 3.2]{An-diagonal}. For the braiding on $r$, we claim that $c( r \ot r ) = r \ot r$. Indeed $r$ is homogeneous of degree $g_1^{12} g_4^6$ and
\begin{align*}
 g_1^{12} g_4^6 \cdot \hu{1j} = (-q_2^2 \omega^2)^6 g_1^{12} \cdot \hu{1j} = (-q_2^2 \omega^2)^6 q_1^{12} \hu{1j} = \hu{1j}, && j\in \I_{2, 3},
\end{align*}
so we have $g_1^{12} g_4^6 \cdot r = r$. Hence $\ku r$ is a braided vector subspace corresponding to an infinite dimensional Nichols algebra, and 
$r \in\toba (\mP (\wtoba))$, a contradiction.  Thus $I = \cJ (\nombre)$ and $\wtoba \simeq \toba(\nombre)$.

\medspace

Now we prove \ref{item:postnichols}. Let us compute $\nombre^* \in \ydg$. Denote by $ ( f_i ) _{ i \in \I_4}$ the basis of $\nombre^*$ dual to $( x_i ) _{ i \in \I_4}$. Then $ \delta (f_i) = g_i^{-1} \ot f_i$, so $c( f_i \ot f_j ) = g_i^{-1} \cdot f_j \ot f_i$. Straightforward computations shows that 
\begin{align}\label{eq:table-dual-action}
&\begin{aligned}
g_i ^ {-1} \cdot f_j &= -f_{2i - j}, & g_i ^ {-1} \cdot f_4 &= q_1 f_4, & \\
g_4 ^ {-1} \cdot f_j &= q_2 f_j, & g_4 ^ {-1} \cdot f_4 &= -\omega^2 f_4,
\end{aligned} & i, j &\in \mathbb{I}_3.
\end{align}
Hence $\nombre \to \nombre ^ *$, $x_i \mapsto f_i$, $i \in \I_4$, is an isomorphism of braided vector spaces.

Let $\cL =\oplus_{n\geq0} \cL^n$ be a finite-dimensional post-Nichols algebra of $\nombre$. By \cite[Lemma 5.5]{AS-adv}, $\cL^*$ is a finite-dimensional pre-Nichols algebra of $\nombre^*\simeq \nombre$. By \ref{item:prenichols}, $\cL^*\simeq \toba(\nombre^*)$, hence $\cL\simeq \toba(\nombre)$.
\epf

\begin{theorem}\label{thm:gen-degree-one}
Let $H$ be a finite-dimensional pointed Hopf algebra over $\Gamma$ with infinitesimal braiding $\nombre$. 
Then $H$ is generated by its group-like and skew-primitive elements.
\end{theorem}
\pf
Working as in \cite[Theorem 5.5]{AS-Annals}, we reduce to prove that the unique finite-dimensional post-Nichols algebra of $\nombre$ is the Nichols algebra, which follows by Theorem \ref{thm:pre-post-Nichols} \ref{item:postnichols}.
\epf

\subsection{The dintinguished pre-Nichols algebra}

Now we introduce a pre-Nichols algebra of $\nombre$ which mimics those given  in \cite{An-distinguished} for braidings of diagonal type.

\begin{definition}
The \emph{distinguished pre-Nichols algebra} $\wtoba(\nombre)$ of $\nombre$ is the quotient of $T(\nombre)$ by the ideal $\cJt(\nombre)$ generated by the elements
\eqref{eq:HV-xi^2}, \eqref{eq:HV-sum-xixj}, \eqref{eq:HV-ui1=u1k} and \eqref{eq:HV-link-x4-zh}.
\end{definition}

\begin{rem}\label{rem:Z-HV}
$\cJt(\nombre)$ is a Hopf ideal by Remark \ref{rem:comultiplication-relations}, and there exists a canonical projection
$\pi:\wtoba(\nombre) \twoheadrightarrow \toba(\nombre)$. Also, \eqref{eq:HV-zh^2}--\eqref{eq:HV-link-xi-u1j} hold in $\wtoba(\nombre)$ by Lemma \ref{lemma:discard-rels}. Let $\cZ(\nombre)$ be the subalgebra of $\wtoba(\nombre)$ generated by 
\begin{align}
z_4 &:=x_4^6 && \text{and} & z_{124134}&:= \hu{12413}^{3}.
\end{align} 
Then $\toba(\nombre)$ is the quotient of $\wtoba(\nombre)$ by the ideal generated by $\cZ(\nombre)$. 
\end{rem}

Next we will prove that $\cZ(\nombre)$ is a normal Hopf subalgebra of $\wtoba(\nombre)$. In order to do so, we need an auxiliar computation.

\begin{lemma} \label{lemma:adx4u1j}
Let $j\in\I_{2, 3}$. The following relations hold in $\wtoba(\nombre)$:
\begin{align} \label{eq:adx4u1j}
(\ad_c x_4) \hu {1j} &= q_2 \omega^2 \z{1} \z{j} - q_2 \z{5-j \,} \z{1},
\\ \label{eq:u1jx4-x4u1j}
\hu{1j} x_4 + q_1^2 \omega ^2 x_4 \hu{1j} &= q_2^{-1} \z{1} \z{j} + q_1 \z{5-j \,} \z{1}.
\end{align}
\end{lemma}
\pf
Fix $j \in \I_{2, 3}$. For \eqref{eq:adx4u1j}, we compute
\begin{align*}
(\ad_c & x_4) \hu{1j} = x_4 \hu{1j} - (g_4 \cdot \hu{ 1 j}) x_4 = x_4 \hu{1j} + q_2 ^2 \omega ^ 2  \hu{ 1 j} x_4 \\
&= x_4 (x_1 \z{j} + q_1 \z{5-j \,} x_1) + q_2^2 \omega^2 (x_1 \z{j} + q_1 \z{5-j \,} x_1) x_4 \\
&= x_4 x_1 \z{j} + q_1 q_2 \z{5-j \,} x_4 x_1 + q_2 \omega^2 x_1 x_4 \z{j} - q_2 \z{5-j \,} ( \z{1} + q_1 x_4 x_1) \\
&= q_2 \omega^2 (x_1 x_4 - q_1 x_4 x_1) \z{j} - q_2 \z{5-j \,} \z{1} \\ %+ (q_1 q_2 - q_1 q_2) \z{5-j \,} x_4 x_1
&= q_2 \omega^2 \z{1} \z{j} - q_2 \z{5-j \,} \z{1}.
\end{align*}
Now \eqref{eq:u1jx4-x4u1j} follows from \eqref{eq:adx4u1j} by multiplying both sides   by 
$q_1^2 \omega ^2$.
\epf

\begin{lemma}\label{lem:Z-HV-str}
\begin{enumerate}[leftmargin=*,label=\rm{(\alph*)}]
\item\label{item:Z-HV-central} $\cZ(\nombre)$ is a normal Hopf subalgebra of $\wtoba(\nombre)$.
\item\label{item:Z-HV-coinv} $\cZ(\nombre) = ^{\co \pi} \wtoba(\nombre) = \wtoba(\nombre)^{\co \pi}$.
\item\label{item:Z-HV-qpol} $\cZ(\nombre)$ is presented by generators $z_4$, $z_{124134}$, which satisfy the relation 
\begin{align*}
z_4z_{124134}=q_2^{72} z_{124134}z_4.
\end{align*} 
The set $\widetilde{B}=\{ z_4^mz_{124134}^n \colon m,n\in\N_0 \}$ is a basis of $\cZ(\nombre)$.
\end{enumerate}
\end{lemma}
\pf
\ref{item:Z-HV-central} We show that the generators of $\cZ(\nombre)$ are vanished by the braided adjoint action of $\wtoba(\nombre)$.

We start with $z_4$. Clearly $(\ad_c x_4)z_4= x_4 x_4^6-(-\omega^2)^6 x_4^6 x_4 = 0$. Let $i\in\I_3$. 
By \cite[Equation (A.8)]{AS-adv} we have
\begin{align*}
(\ad_c x_4)^6 x_i &= \sum_{k=0}^{6} (-1)^k \binom{6}{k}_{-\omega^2} (-\omega^2)^{k(k-1)/2} q_2^{k} \, x_4^{6-k}x_ix_4^{k}  = x_4^{6}x_i - q_2^{6} x_ix_4^{6} \\
&=-q_2^{6}(x_i x_4^6-q_1^{6} x_4^{6} x_i )=-q_2^{6} (\ad_c x_i)z_4.
\end{align*}
On the other hand,  $(\ad_c x_4)^2 x_i = 0$ by \eqref{eq:HV-link-x4-zh}, hence $(\ad_c x_i)z_4 = 0$.

Now turn to $z_{124134}$. By Lemma \ref{lemma:discard-rels}, we have $\hu{1 j}^2 = 0$ for $j \in \I_{2, 3}$, so 
\begin{align*}
z_{124134} &= (\hu{12} \hu{13} + \omega^2 \hu{13} \hu{12})^3 
\\
&= \hu{12} \hu{13} \hu{12} \hu{13} \hu{12} \hu{13} + \hu{13} \hu{12}\hu{13} \hu{12}\hu{13} \hu{12}.
\end{align*}
Let $i\in\I_3$. By \eqref{eq:HV-link-xi-u1j}, $(\ad_c x_i)\hu{1 j} = 0$ for $j\in\I_{2,3}$. Since $\ad_c x_i$ is a braided derivation, it follows that $(\ad_c x_i)z_{124134}=0$.
%Using \eqref{eq:uij-trenza} and \eqref{eq:HV-link-xi-u1j} we have
%\begin{align*}
%(\ad_c x_i)z_{124134} = x_i z_{124134} - (g_i \cdot z_{124134}) x_i 
%= x_i z_{124134} - q_1^6 z_{124134} x_i =0.
%\end{align*}

By Lemma \ref{lemma:discard-rels}, \eqref{eq:HV-link-u1j-zh} holds in $\wtoba(\nombre)$ since this relation follows from \eqref{eq:HV-zh^2} and \eqref{eq:HV-uij=u1k}. By repeated applications of \eqref{eq:u1jx4-x4u1j} and \eqref{eq:HV-link-u1j-zh} we get
\begin{align*}
& z_{124134} x_4= \hu{12} \hu{13} \hu{12} \hu{13} \hu{12} (-q_1^2 \omega^2 x_4 \hu{13} + q_2 ^{-1} \z{1} \z{3} + q_1 \z{2} \z{1}) + \\
& \qquad + \hu{13} \hu{12}\hu{13} \hu{12}\hu{13} (-q_1^2 \omega^2 x_4 \hu{12} + q_2 ^{-1} \z{1} \z{2} + q_1 \z{3} \z{1})  \\
& \quad = -q_1^2 \omega^2 \hu{12} \hu{13} \hu{12} \hu{13} (-q_1^2 \omega^2 x_4 \hu{12} + q_2 ^{-1} \z{1} \z{2} + q_1 \z{3} \z{1}) \hu{13}
\\
& \qquad + q_2^{-1} q_1^{10} \z{2} \z{1}\hu{12} \hu{13} \hu{12} \hu{13} \hu{12} + q_1^{11} \z{3} \z{2} \hu{12} \hu{13} \hu{12} \hu{13} \hu{12} 
\\ 
& \qquad -q_1^2 \omega^2 \hu{13} \hu{12} \hu{13} \hu{12}(-q_1^2 \omega^2 x_4 \hu{13} + q_2 ^{-1} \z{1} \z{3} + q_1 \z{2} \z{1}) \hu{12}
\\
& \qquad + q_2^{-1} q_1^{10} \z{3} \z{1}\hu{13} \hu{12} \hu{13} \hu{12}\hu{13} + q_1^{11} \z{2} \z{3} \hu{13} \hu{12} \hu{13} \hu{12}\hu{13} 
\\
& \quad = \dots =(-q_1^2 \omega^2 )^6 x_4 z_{124134} + 
\\
& \qquad +q_1^{11}( \z{1} \z{2} -(\omega +\omega^2) \z{3} \z{1} +  \z{2} \z{3})  \hu{13} \hu{12} \hu{13} \hu{12}\hu{13} 
\\
& \qquad +q_1^{11}(\z{1} \z{3} -(\omega+\omega^2) \z{2} \z{1} + \z{3} \z{2}) \hu{12} \hu{13} \hu{12}\hu{13} \hu{12}.
\end{align*}
By \eqref{eq:HV-sum-zizj} we have that
\begin{align}\label{eq:adx4-z}
(\ad_c x_4) z_{124134} &= %x_4 z_{124134} - (g_4 \cdot z_{124134}) x_4 
x_4 z_{124134} - q_2^{12}z_{124134} x_4 = 0.
\end{align}
Hence $\cZ(\nombre)$ is a normal subalgebra. Also $\cZ(\nombre)$ is a Hopf subalgebra since $z_4, z_{124134} \in \mP(\wtoba(\nombre))$ by Remark \ref{rem:comultiplication-relations} and \cite[Lemma 3.2]{An-diagonal}

\ref{item:Z-HV-coinv} This fact follows from \cite[Proposition 3.6 (c)]{A+}.

\ref{item:Z-HV-qpol} Note that $z_4, z_{124134} \neq 0$ since \eqref{eq:HV-xi^2}-\eqref{eq:HV-link-x4-zh} minimally generate $\cJ(\nombre)$.
Also, $z_4$ and $z_{124134}$ $q$-commute by \eqref{eq:adx4-z}, so $\cZ(\nombre)$ is spanned by $\widetilde{B}$. 

Let $K$ be the subalgebra of $\wtoba(\nombre)\# \Bbbk\Gamma$ generated by $z_4, z_{124134}$ and $\Gamma$: $K$ is a Hopf subalgebra, which is a pointed Hopf algebra since $z_4$ is $(1,g_4^6)$-primitive and
$z_{124134}$ is $(1,g_1^{12}g_4^6)$-primitive. 
As $z_4, z_{124134}$ are linearly independent, the infinitesimal braiding contains the braided vector space generated by them, which is of diagonal type with matrix $\begin{pmatrix}
1 & q_1^{72} \\ q_2^{72} & 1 
\end{pmatrix}$ so the set $\{z_4^m z_{124134}^n \gamma \colon m,n\in\N_0, \gamma\in\Gamma \}$ is linearly independent. Thus $\widetilde{B}$ is linearly independent.
\epf

\begin{pro}
\begin{enumerate}[leftmargin=*,label=\rm{(\alph*)}]
\item\label{item:prenichols-extension} We have an extension of braided Hopf algebras:
\begin{align*}
\Bbbk \to \cZ(\nombre) \hookrightarrow \wtoba(\nombre) \twoheadrightarrow \toba(\nombre) \to \Bbbk.
\end{align*}

\item\label{item:prenichols-GKdim} $\GK \wtoba(\nombre) = 2$.

\item\label{item:prenichols-PBW} The following set is a basis of $\wtoba(\nombre)$:
\begin{align}\label{dist:pbw}
\begin{aligned}
x_4^{n_4}
\left\lbrace 
\begin{aligned}
\z{1}^{n_{14}}\z{2}^{n_{24}}\z{3}^{n_{34}}\\
\z{1}^{n_{14}}\z{2}\z{1}\z{3}^{n_{34}}
\end{aligned}
\right\rbrace 
\hu{12}^{n_{124}}\hu{12413}^{n_{124134}}\hu{13}^{n_{134}}
\left\lbrace 
\begin{aligned}
x_1^{n_1}x_2^{n_2}x_3^{n_3}\\
x_1^{n_1}x_2x_1x_3^{n_3}
\end{aligned}
\right\rbrace,
\\
n_{14}, n_{24}, n_{34}, n_{124}, n_{134}, n_1, n_2, n_3 \in\I_{0,1}, \, n_4, n_{124134}\in\N_0.
\end{aligned}
\end{align}
\end{enumerate}
\end{pro}
\pf
\ref{item:prenichols-extension} This fact follows by Lemma \ref{lem:Z-HV-str} and Remark \ref{rem:Z-HV}, cf. \cite[\S 2.5]{AN}.

\smallbreak
\ref{item:prenichols-GKdim} 
By \cite[Proposition 3.6 (d)]{A+} there exists a right $\cZ(\nombre)$-linear isomorphism $\toba(\nombre)\ot \cZ(\nombre) \simeq \wtoba(\nombre)$, hence $\wtoba(\nombre)$ is finitely generated as right $\cZ(\nombre)$-module. By 
\cite [Proposition 5.5]{KL}, 
\begin{align*}
\GK \wtoba(\nombre) = \GK  \cZ(\nombre) = 2.
\end{align*}

\smallbreak
\ref{item:prenichols-PBW} First, $\wtoba(\nombre)$ is spanned by the set \eqref{dist:pbw}: this follows using the defining relations and arguing as in Lemma \ref{lemma:HV-presentation}. 

Let $\{b_i\}_{i\in\I_{10368}}$ be an enumeration of \eqref{D3-1:pbw}. Let $\widetilde{b}_i$ be the corresponding element viewed in $\wtoba(\nombre)$; thus $\pi(\widetilde{b}_i)=b_i$. Since the elements of $\cZ(\nombre)$ $q$-commute with all the elements of $\wtoba(\nombre)$, the set \eqref{dist:pbw} coincides with 
$$ B=\{z_4^{m} z_{124134}^{n} \widetilde{b}_i : m,n \in\N_0, i\in\I_{10368} \} $$ 
up to non-zero scalars. Hence it suffices to prove that $B$ is linearly independent. Suppose that $0=\sum_{m,n \in\N_0} \sum_{i\in\I_{10368}} a_{mni} z_4^{m} z_{124134}^{n} \widetilde{b}_i$, where not all the scalars $a_{mni}$'s are zero. Let $i_0\in \I_{10368}$ be such that $a_{mni_0} \neq 0$ for some $m,n\in\N_0$ and $b_i$ is of maximal degree $N$. Let $f\in\toba(\nombre)^*$, $f(b_i)=\delta_{ii_0}$, $i\in \I_{10368}$. Since $\Delta$ and $\pi$ are $\N_0$-graded, 
$(\id \ot f)(\id \ot\pi)\Delta(\widetilde{b}_i)=0$ for all $\widetilde{b}_i$ of degree less than $N$; for those $\widetilde{b}_i$ of degree $N$,
$(\id \ot f)(\id \ot\pi)\Delta(\widetilde{b}_i)=\delta_{ii_0}1$ since 
$\Delta(\widetilde{b}_i)\in 1\ot \widetilde{b}_i+ \sum_{j>0} \wtoba(\nombre)^j \ot \wtoba(\nombre)^{N-j}$.
Hence
\begin{align*}
0 &= (\id \ot f) (\id \ot\pi)\Delta \Big(\sum_{m,n \in\N_0} \sum_{i\in\I_{10368}} a_{mni} z_4^{m} z_{124134}^{n} \widetilde{b}_i\Big)
\\
&= \sum_{m,n \in\N_0} \sum_{i\in\I_{10368}} a_{mni} (\id \ot f) (\id \ot\pi) \Big(\Big(\sum_{j\in\I_{0,m}} \binom{m}{j} z_4^j \ot z_4^{m-j}\Big)
\\ & \qquad \times \Big(\sum_{k\in\I_{0,n}} \binom{n}{k} z_{124134}^k \ot z_{124134}^{n-k}\Big)\Delta(\widetilde{b}_i)\Big)
\\
&= \sum_{m,n \in\N_0} \sum_{i\in\I_{10368}} a_{mni} (\id \ot f)  
\Big((z_4^m z_{124134}^n \ot 1)
(\id \ot\pi)\Delta(\widetilde{b}_i)\Big)
\\
&= \sum_{m,n \in\N_0} \sum_{i\in\I_{10368}} a_{mni} (z_4^m z_{124134}^n \ot 1) (\id \ot f)(\id \ot\pi)\Delta(\widetilde{b}_i)
\\
&= \sum_{m,n \in\N_0} a_{mni_0} z_4^m z_{124134}^n,
\end{align*}
which is a contradiction by Lemma \ref{lem:Z-HV-str} \ref{item:Z-HV-qpol}.
\epf

\section{Liftings of $\nombre$} \label{sec:liftings}
We fix a group $\Gamma$ and a principal realization of $\nombre$ as in \S \ref{subsec:Gamma3}.
In order to compute the liftings of $\toba ( \nombre )$ over $\Gamma$ we follow the strategy developed in \cite{A+,AAG}.
The procedure gives rise to a family of liftings which are cocycle deformations of $\toba(\nombre) \# \Bbbk\Gamma$. 
Moreover, it provides a criterion to check if we have an exhaustive list of liftings of $\toba(\nombre)$ over $\Gamma$.\\

The starting point of the strategy is a suitable chosen chain of subsequent quotients of pre-Nichols algebras of $\nombre$:
\begin{align*}
T(\nombre) =: \toba_0 \twoheadrightarrow \toba_1 \twoheadrightarrow \cdots \twoheadrightarrow \toba_{l + 1} = \toba (\nombre).
\end{align*}
After bosonization with $\Bbbk\Gamma$ we obtain graded Hopf algebras $\mH _k := \toba_k \# \Bbbk\Gamma$, $k \in \I_{l+1}$, related by a corresponding chain of quotients.
If $\cA  \in \Cleft (\mH_k)$, then the left Schauenburg Hopf algebra $L(\cA,\mH_k)$ is (isomorphic to) a quotient of $\mT:=\mH_0 =T(\nombre)\#\Bbbk\Gamma$ by \cite[Proposition 5.10]{A+}; this yields a filtration $\mathfrak F$ on $L(\cA , \mH_k)$, induced by the filtration coming from the $\N_0$-graduation of $T(\nombre)$. We denote by $\gr_{\mathfrak F}$ the associated graded object. Let
\begin{align*}
\Cleft'(\mH _k) &: = \{ \cA \in \Cleft (\mH_k) \colon  \gr_{\mathfrak F} L(\cA,\mH_k) \simeq \mH_k \}, & k& \in \I_{l+1}.
\end{align*}
Recursively, the strategy begins with $\Cleft'(\mT) = \{\mT\}$ and compute the set $\Cleft'(\mH_{k+1})$ from $\Cleft'(\mH_k)$ using the ideas in \cite{Gu}.

Since $\toba(\nombre)$ is coradically graded, if $L$ is a lifting of $\nombre$ then $\mathfrak F $ above coincides with the coradical filtration and we have 
\begin{align*}
\Cleft'(\mH)=\{\cA\in\Cleft(\mH)\colon\gr L(\cA,\mH) \simeq \mH  \},
\end{align*}
so the last step of the recursion leaves us with a family of liftings that are cocycle deformations of $\toba(\nombre)$.

\medbreak

We mention particular features of the strategy that will clarify the upcoming computations. Fix $k\in\I_{l+1}$.
\begin{itemize}[leftmargin=*]
\item All cleft objects of $\mH_{k+1}$ arise as quotients of cleft objects of $\mH_{k}$ \cite{Gu}.
\item If $ \cA\in \Cleft'(\mH_k )$, there is an algebra map $\mT 
\twoheadrightarrow \cA$. The $\mH_k$-colinear section $\gamma_k \colon \mH_k \to \cA $ restricts to an
algebra map $(\gamma_k) _{|\Bbbk\Gamma} \in \Alg (\Bbbk\Gamma, \cA)$.
\item Each $\cA\in\Cleft'(\mH_k )$ can be obtained as $\cA=\cE\#\Bbbk\Gamma$ where $\cE$ is a $\Bbbk\Gamma$-module algebra \cite[Proposition 5.8]{A+}.
Moreover, the section $\gamma_k$ restricts to a braided $\toba_k$-comodule isomorphism $\gamma_k \colon \toba_k \to \cE$
\cite {AnG2}.
\item As algebras, each $\cE$ is a quotient of some $\cE'$ of the previous step.
\item $\cL_0 = \mT$ and, if $\cA\twoheadrightarrow \cA'$, then $L (\cA,\mH_{k }) \twoheadrightarrow L (\cA', \mH_{ k + 1})  $.
\end{itemize}

\subsection{Stratification}
The first task proposed by the strategy is to obtain a convenient stratification $\cG= \cG_0 \sqcup \cG_1 \sqcup \cdots  \sqcup \cG_l$ 
of the minimal set of generators of the ideal $\cJ (\nombre)$ found in Theorem \ref{thm:minimal-presentation}. 
We need that the elements of $\cG_k$ are primitive in the braided Hopf algebra $\toba_k := T ( \nombre ) / \langle \cup_{j=0}^{k-1} \cG_j \rangle$, $k  \in \I_{l+1}$. In our setting, we also require that the vector space spanned by
$\cG_k$ is a Yetter-Drinfeld submodule of $T(\nombre)$ over $\Bbbk\Gamma$.

By Remark \ref{rem:comultiplication-relations} we may consider the following stratification:
\begin{align*}
\cG_0=& \{x_4 \z{h} - q_2 \z{h} x_4 \colon  h \in \I_3 \}  ,& \cG_1&= \{ x_h^2 \colon h \in \I_3 \} ,
\\
\cG_2 =& \{ x_1 x_i + x_{5-i} x_1 + x_i x_{5-i} \colon i \in \I_{2, 3} \}, & \cG_3&= \{ \hu{i1} - \omega \hu{1 5-i \,} \colon i\in \I_{2,3} \},
\\
\cG_4= & \{ x_4^6, \, (\hu{12}\hu{13}+\omega^2\hu{13}\hu{12})^3 \} .
\end{align*}

\subsubsection{Realization of the strata}

To describe the liftings we need to determine the braided vector space structure of each step of the stratification. To this end, it is enough to describe the action of $g_i\in\Gamma$, $i\in\I_4$.

\smallbreak

Recall (cf. \S\ref{subsec:Gamma3}) that we have elements $g_i \in \Gamma$, 1-cocycles $\chi_i:\Gamma\to\ku$, $i \in \I_4$, and an action $\triangleright \colon \Gamma \times \I_4 \to \I_4$  such that
\begin{align*}
x_h \in T(\nombre)_{g_h},&&g \cdot x_h &= \chi_h (g) x_{g \triangleright h}, && g \in \Gamma, \, h \in \I_4.
\end{align*}
We summarize the structure of the submodules $\ku \cG_k \subset T(\nombre) \in \ydg$, $k \in \I_{0,4}$, for later reference.\\
Let $\widetilde{\triangleright} \colon\I_4 \times \I_{2,3} \to \I_{2,3}$ such that 
$j \, \widetilde{\triangleright} \, i =\begin{cases} 5-i &\text { if } j \in \I_3,\\ i & \text{ if } j=4.\end{cases}$

\noindent $\heartsuit$ $\cG_0$: structure determined by $x_4 \z{h} - q_2 \z{h} x_4\in T(\nombre)_{g_hg_4^2}$, $h\in\I_3$,
\begin{align*}
g\cdot(x_4 \z{h} - q_2 \z{h} x_4) &=  \chi_h \chi_4^2(g) (x_4 \z{g \triangleright h} - q_2 \z{g \triangleright h} x_4),
&& g \in \Gamma, \,h \in \I_3.
\end{align*}

\medbreak
\noindent $\heartsuit$ $\cG_1$: Here $x_h^2\in T(\nombre)_{g_h^2}$ and the action satisfies 
\begin{align*}
g\cdot x_h^2 &=  \chi_h^2 (g)x_{g \triangleright h}^2, & g \in \Gamma, & \,h \in \I_3.
\end{align*}

\medbreak
\noindent $\heartsuit$ $\cG_2$. Put $r_i:=x_1 x_i + x_{5-i} x_1 + x_i x_{5-i}$, $i \in \I_{2,3}$. In this case $r_i \in T(\nombre)_{g_1g_i}$ and the action satisfies
\begin{align*}
g_j \cdot r_i = \chi_1\chi_i(g_j)r_{j \,\widetilde{\triangleright} \, i}, && j\in \I_4, \, i \in \I_{2, 3}.
\end{align*}

\medbreak
\noindent $\heartsuit$ $\cG_3$. Let $p_i:=\hu{i1} - \omega \hu{1 5-i \,}$. Then $p_i\in T(\nombre)_{g_1g_{5-i}g_4}$, $i \in \I_{2,3}$. By \eqref{eq:action-ui1=u1k}, there are maps $\eta_i \colon \I_4 \to \ku^{\times}$, $i \in \I_{2,3}$ such that
\begin{align}
g_j \cdot p_i= \eta_i(j) p_{j \,\widetilde{\triangleright} \, i }, && j\in \I_4, \, i \in \I_{2, 3}.
\end{align}

\medbreak
\noindent $\heartsuit$ $\cG_4$.  Here $x_4^6 \in T(\nombre)_{g_4^6}^{\chi_4^6}$, and $(\hu{12}\hu{13}+\omega^2\hu{13}\hu{12})^3 \in T(\nombre)_{g_1^{12}g_4^6}^{\chi_1^{12}\chi_4^6}$.

\medspace

For $r\in\cG$, we denote by $g_r$ the element of $\Gamma$ such that $r\in T(\nombre)_{g_r}$.

\begin{lemma}\label{lem:gr-neq-gi}
If $r\in\cG$, then $g_r \neq g_i$ for all $i \in \I_4$.
\end{lemma}
\pf
Suppose first $r \in \cG_0$, so $g_r = g_hg_4^2$ for some $h \in \I_3$. Since $g_4$ is central and $g_h$ is not, we have $g_hg_4^2\neq g_4$. Also $g_hg_4^2\neq g_h$, because $g_4^2 \cdot x_4 = \omega x_4$ and $\omega \neq 1$. If $g_hg_4^2=g_i$ for some $i \in \I_3$, $i\neq h$, then $-x_i=g_i\cdot x_i=g_hg_4^2\cdot x_i=-q_2^2 x_{2h-i}$, a contradiction. 

\medspace
Next, let $r \in \cG_1$, so $g_r=g_h^2$ for some $h\in \I_3$. Since $g_h^2=g_i^2$ for all $i \in \I_3$, it follows $g_h^2\neq g_i$, $i \in \I_3$. If $g_h^2=g_4$ then 
\begin{align*}
x_h=g_h^2\cdot x_h = g_4\cdot x_h = q_2x_h, && -\omega^2x_4=g_4\cdot x_4 = g_h^2\cdot x_4=q_1^2x_4,
\end{align*}
hence $\omega=(q_1q_2)^{10}=-\omega$, contradicting $\omega \in \Gp_3$.

\medspace
Suppose $r \in \cG_2$, so $g_r=g_1g_i$ for some $i \in \I_{2,3}$. Clearly $g_1g_i\neq g_1, g_i$; since $g_{5-i}g_1=g_1g_i$, it also follows $g_1g_i\neq g_{5-i}$. If $g_1g_i = g_4$ then we have $x_i=g_1g_i \cdot x_1 = g_4\cdot x_1=q_2x_1$, a contradiction.

\medspace
Assume now $r\in \cG_3$ and let $i\in\I_{2,3}$ such that $g_r=g_1g_{5-i}g_4$. Since $g_1g_{5-i}=g_ig_1=g_{5-i}g_i$, we have $g_1g_{5-i}g_4\neq g_1, g_i, g_{5-i}$, because $g_{5-i}, g_1, g_i$, respectively, are non-central. We also have $g_1g_{5-i}g_4\neq g_4$, since $g_1g_{5-i}$ acts non-trivially on $x_1$.

\medspace
Turn now to $\cG_4$. We have $g_4^6\neq g_h$ for $h\in\I_3$ because $g_h$ is non-central; also $g_4^4\neq g_4$, because $g_4^5\cdot x_4=-\omega x_4$ and $-\omega\neq1$. Finally, assume $g_r=g_1^{12}g_4^6$. If $g_1^{12}g_4^6 = g_h$ for some $h\in\I_3$, then for any $j \in \I_3$ different from $h$ we have  $-x_{2h-j}=g_h\cdot x_j=g_1^{12}g_4^6 \cdot x_j=g_j^{12}g_4^6 \cdot x_j \in \ku x_j$, a contradiction. If $g_1^{12}g_4^6 = g_4$, we compute
\begin{align*}
x_1= g_1^{12}g_4^5\cdot x_1=q_2^5, && x_4=g_1^{12}g_4^5\cdot x_4=-q_1^{12}\omega x_4,
\end{align*}
so $q_2^5=1$, $q_1^{12}=-\omega^2$. Then $1=(q_1q_2)^{60}=(-\omega^2)^5= -\omega$, a contradiction.
\epf

\subsection{Computing cleft objects}

The second task of the strategy is the introduction of a suitable family of cleft objects of $\mH$. To this end we define a family of $\Bbbk\Gamma$-module algebras such that, after bosonization with $\Bbbk\Gamma$, give the desired cleft extensions.

\smallbreak

The set of deformation parameters $\cR_{\nombre}$ consists of 4-uples of scalars $\bsl=(\lambda_i)_{i\in\I_4} \in\Bbbk^4$ such that 
\begin{align}\label{eq:lambda-conditions}
\begin{aligned}
\lambda_1 &= 0 \text { if either } \chi_i^2 \neq \varepsilon \text { or } g_i^2 = 1 \text { for some } i\in\I_3, 
\\
\lambda_2 &= 0 \text { if either } \chi_i \chi_j  \neq \varepsilon \text { or } g_i g_j = 1 \text { for some } i\neq j\in\I_3,
\\
\lambda_3 &= 0 \text { if either } \chi_4^6 \neq \varepsilon \text { or } g_4^6=1,
\\
\lambda_4 &= 0 \text { if either } \chi_1^{12} \chi_4^6 \neq \varepsilon \text{ or } g_1^{12} g_4^6=1.
\end{aligned}
\end{align}

Let $\bsl\in \cR_{\nombre}$. We define $\cE_0(\bsl) =\toba_0=T(\nombre)$, $\cE_1(\bsl) =\toba_1$, but we change the labels of the generators to $(y_i)_{i\in\I_4}$ in order to differentiate with generators $(x_i)_{i\in\I_4}$ of the pre-Nichols algebras $\toba_k$. Let
\begin{align*}
\cE_{i+1}(\bsl) & := \cE_{i}(\bsl) / \left\langle r - \lambda_i \colon r\in\cG_i\right\rangle, & &i\in\I_{2},
\\
\cE_4(\bsl) &:= \cE_3(\bsl)/ \left\langle r - \lambda_2 y_4 \colon r\in\cG_3\right\rangle.
\end{align*} 

\begin{remark}
Each $\cE_{i}(\bsl)$ is a $\Bbbk\Gamma$-module algebra since the ideal is stable by the $\Gamma$-action by \eqref{eq:lambda-conditions}. Thus we may introduce $\cA_i(\bsl):=\cE_i(\bsl) \# \Bbbk\Gamma$.
\end{remark}

\medbreak

\begin{lemma}
Let $k\in\I_{4}$. Then $\cE_k(\bsl)\neq 0$ and each $\cA_k(\bsl)$ is a $\mH_k$-cleft object. There exists an $\mH_k$-colinear section $\gamma_k \colon \mH_k \to \cA_k$ which restricts to an
algebra map $(\gamma_k) _{|\Bbbk\Gamma} \in \Alg (\Bbbk\Gamma, \cA_k)$.
\end{lemma}
\pf
Fix $\bsl\in\cR(\nombre)$; we prove the claim recursively on $k$. For the sake of simplicity of the notation, we call $\cE_k=\cE_k(\bsl)$, $\cA_k=\cA_k(\bsl)$.
For $k=1$ the claim is clear since $\cE_1=\toba_1$, $\cA_1=\mH_1$ so we take $\gamma_1=\id_{\mH_1}$.

\smallbreak

\noindent $\heartsuit$ 
For $k=2$, we notice that $\cE_2 \neq 0$ (and a fortiori $\cA_2\neq 0$) since 
there exists an algebra map $\cE_2 \twoheadrightarrow \cE(\lambda_1,\lambda_2)$, where   $\cE(\lambda_1,\lambda_2)$ is the corresponding non-trivial algebra given in
\cite[Proposition 7.2]{GV2}: the map identifies the generators $y_i$ for $i\in\I_3$ and annihilates $y_4$. As 
\begin{align*}
g_3(y_1^2-\lambda_1)g_3^{-1}&=y_2^2-\lambda_1, & g_2(y_1^2-\lambda_1)g_2^{-1}&=y_3^2-\lambda_1,
\end{align*}
we have that $\langle y_i^2 -\lambda_1: i\in\I_3 \rangle=\langle y_1^2 -\lambda_1 \rangle$. Using \cite[Proposition 5.8]{A+} and this equality of ideals, we have that
\begin{align*}
\cA_2 &\simeq \big(\cE_1\# \Bbbk\Gamma\big)/ \langle y_i^2 -\lambda_1: i\in\I_3 \rangle
= \cA_1/ \langle y_1^2 -\lambda_1 \rangle.
\end{align*}

Let $Y_1'$ be the subalgebra of $\mH_1$ generated by 
$x_1^2$. Then $Y_1'$ is isomorphic to a polynomial ring in one variable since $x_1^2 \in \mP(\toba_1)_{g_1^2}-0$ and $g_1^2\cdot x_1^2=x_1^2$. As in \cite{AAG} we set $Y_1=\Ss(Y_1')$. Notice that
\begin{align*}
\mH_1/ \langle Y_1^+ \rangle = \mH_1/ \langle x_1^2 \rangle = \mH_1/ \langle x_i^2  \colon i\in \I_3 \rangle \simeq \mH_2. 
\end{align*}
Since $Y_1$ is also a polynomial algebra generated by $x_1^2g_1^{-2}$, there exists an algebra map $\phi:Y_1\to\cA_1$ such that
$\phi(x_1^2g_1^{-2})=y_1^2g_1^{-2}-\lambda_1g_1^{-2}$, which is $\mH_1$-colinear. We notice that 
\begin{align*}
\cA_1/ \langle \phi(Y_1^+) \rangle &= 
\cA_1/ \langle y_1^2-\lambda_1 \rangle \simeq  \cA_2.
\end{align*}
Hence $\cA_2$ is a $\mH_2$-cleft object by \cite[Theorem 8]{Gu}. The claim about the section $\gamma_2$ follows from \cite[Proposition 5.8]{A+}.

\smallbreak

\noindent $\heartsuit$ 
For $k=3$, $\cE_3 \neq 0$ (and a fortiori $\cA_3\neq 0$) since 
the algebra map $\cE_2 \twoheadrightarrow \cE(\lambda_1,\lambda_2)$ descends to an algebra map $\cE_3 \twoheadrightarrow \cE(\lambda_1,\lambda_2)$. 

Let $\rel_i= y_1y_i + y_{5-i}y_1+y_iy_{5-i}$, $i\in\I_{2, 3}$. As $g_1(\rel_2-\lambda_2)g_1^{-1}=\rel_3-\lambda_2$,
we have that $\langle \rel_i -\lambda_2: i\in\I_{2,3} \rangle=\langle \rel_2 -\lambda_2 \rangle$. Hence
\begin{align*}
\cA_3 &\simeq \big(\cE_2\# \Bbbk\Gamma\big)/ \langle \rel_i -\lambda_2: i\in\I_{2,3} \rangle
= \cA_2/ \langle \rel_2 -\lambda_2 \rangle.
\end{align*}

Let $Y_2'$ be the subalgebra of $\mH_2$ generated by 
$r_2=x_1x_2+x_3x_1+x_2x_3$. Then $Y_2'$ is isomorphic to a polynomial ring in one variable since $r_2 \in \mP(\toba_2)_{g_1g_2}-0$ and $g_1g_2\cdot r_2=r_2$. As in \cite{AAG} we set $Y_2=\Ss(Y_2')$. Notice that
$\mH_2/ \langle Y_2^+ \rangle \simeq \mH_3$.
Since $Y_2$ is also a polynomial algebra generated by $r_2g_2^{-1}g_1^{-1}$, there exists an algebra map $\phi:Y_2\to\cA_2$ such that
$\phi(r_2g_2^{-1}g_1^{-1})=\rel_2g_2^{-1}g_1^{-1}-\lambda_2g_2^{-1}g_1^{-1}$, which is $\mH_2$-colinear. 
Hence $\cA_3$ is a $\mH_3$-cleft object by \cite[Theorem 8]{Gu}, since 
$\cA_2/ \langle \phi(Y_2^+) \rangle = 
\cA_2/ \langle \rel_2-\lambda_2 \rangle \simeq  \cA_3$. The claim about the section $\gamma_3$ again follows from \cite[Proposition 5.8]{A+}.

\smallbreak

\noindent $\heartsuit$ 
For $k=4$, we check that $\cE_4 \neq 0$ (and a fortiori $\cA_4\neq 0$) using \texttt{GAP}. 
Let $\relp_i= y_{5-i \, 14}-\omega y_{1i4}-\lambda_2y_4$, $i\in\I_{2, 3}$. As $g_1\relp_2g_1^{-1}=q_1\relp_3$,
we have that $\langle \relp_i: i\in\I_{2,3} \rangle=\langle \relp_2 \rangle$. Hence
$\cA_4 \simeq \cA_3/ \langle \relp_2 \rangle$.

Let $Y_3'$ be the subalgebra of $\mH_3$ generated by 
$p_2=\hu{31}-\omega\hu{12}$. Using \texttt{GAP} we check that 
$p_2^6\neq 0$. Thus, by \cite[Lemma 5.13]{A+}, $Y_3'$ is isomorphic to a polynomial ring in one variable since $p_2 \in \mP(\toba_3)_{g_1g_2g_4}-0$ and $g_1g_2g_4\cdot p_2=-\omega^2 p_2$. As in \cite{AAG} we set $Y_3=\Ss(Y_3')$. Notice that
$\mH_3/ \langle Y_3^+ \rangle \simeq \mH_4$.
Since $Y_3$ is also a polynomial algebra generated by $p_2g_4^{-1}g_2^{-1}g_1^{-1}$, there exists an algebra map $\phi:Y_3\to\cA_3$ such that
$\phi(p_2g_4^{-1}g_2^{-1}g_1^{-1})=\relp_2g_4^{-1}g_2^{-1}g_1^{-1}$, which is $\mH_3$-colinear. 
Hence $\cA_4$ is a $\mH_4$-cleft object by \cite[Theorem 8]{Gu}, since 
$\cA_3/ \langle \phi(Y_3^+) \rangle = 
\cA_3/ \langle \relp_2 \rangle \simeq  \cA_4$. The claim about the section $\gamma_4$ again follows from \cite[Proposition 5.8]{A+}.
\epf

Next we define
\begin{align*}
\cE(\bsl) &:= \cE_4(\bsl)/ \left\langle y_4^{6}-\lambda_3, \,  \gamma_4 \Big( (\hu{12}\hu{13}+\omega^2\hu{13}\hu{12})^3\Big) - \lambda_4\right\rangle.
\end{align*} 
Again, $\cE(\bsl)$ is a $\Bbbk\Gamma$-module algebra, so we may consider $\cA(\bsl):=\cE(\bsl) \# \Bbbk\Gamma$.

\begin{lemma}\label{lem:A-cleft}
$\cA(\bsl)$ is a $\mH$-cleft object. There exists an $\mH$-colinear section $\gamma \colon \mH \to \cA$ which restricts to an
algebra map $(\gamma) _{|\Bbbk\Gamma} \in \Alg (\Bbbk\Gamma, \cA)$.
\end{lemma}
\pf
Notice that $\toba_4=\wtoba(\nombre)$. We keep the notation from Remark \ref{rem:Z-HV}: $z_4 = x_4^6$, $z_{124134}= (\hu{12}\hu{13}+\omega^2\hu{13}\hu{12})^3$. Let $\toba_5 = \toba_4/\langle z_4\rangle$, hence $\toba(\nombre)=\toba_5 / \langle z_{124134} \rangle$. We denote the canonical projections as follows:
\begin{align*}
\pi_1 & \colon \toba_4 \twoheadrightarrow \toba_5, &
\pi_2 & \colon \toba_5 \twoheadrightarrow \toba(\nombre), &
\text{hence }&\pi=\pi_2\circ\pi_1.
\end{align*}
We proceed in two steps, one for each relation.

\smallbreak

\noindent $\heartsuit$ 
By Lemma \ref{lem:Z-HV-str} and \cite[Proposition 3.6 (c)]{A+}, $\Bbbk [z_4] =  \toba_4^{\co \pi_1}$, so $X'_4:= \mH_4^{\co \pi_1\# \id} =\Bbbk [z_4]\# 1$; set $X_4=\Ss(X'_4)$. There exists an algebra map $\phi \colon X_4 \to \cA_4$, $\phi(z_4g_4^{-6})=y_4^6g_4^{-6}-\lambda_3g_4^{-6}$, which is $\mH_4$-colinear. Working as in \cite[Theorems 4.7 \& 5.15]{AAG} we check that this map is also $\mH_4$-linear. Indeed, $y_iy_4^6 = q_1^6 y_4^6 y_i$ for all $i\in\I_3$ since \eqref{eq:HV-link-x4-zh} is not deformed; hence, by \eqref{eq:lambda-conditions},
\begin{align*}
\phi(g_i \cdot z_4g_4^{-6}) &= q_1^6 (y_4^6g_4^{-6}-\lambda_3g_4^{-6}) 
= g_i \cdot \phi(z_4g_4^{-6}),
\\
\phi(x_i \cdot z_4g_4^{-6}) &=0 = y_i(y_4^6g_4^{-6}-\lambda_3g_4^{-6}) - q_1^6 (y_4^6g_4^{-6}-\lambda_3g_4^{-6})y_i 
\\
&=x_i \cdot \phi(z_4g_4^{-6}),
\\
\phi(g_4 \cdot z_4g_4^{-6})&= y_4^6g_4^{-6}-\lambda_3g_4^{-6} = g_4 \cdot \phi(z_4g_4^{-6}),
\\
\phi(x_4 \cdot z_4g_4^{-6})&= 0 = x_4 \cdot \phi(x_4g_4^{-6}).
\end{align*}

Let $\mH_5:=\toba_5\# \Bbbk\Gamma$, $\cE_5:=\cE_4/ \langle y_4^6-\lambda_3 \rangle$, $\cA_5:= \cE_5\# \Bbbk \Gamma$.
By \cite[Theorem 4]{Gu} $\cA_5$ is a $\mH_5$-cleft object.

\smallbreak

\noindent $\heartsuit$ 
By Lemma \ref{lem:Z-HV-str} and \cite[Proposition 3.6 (c)]{A+}, $\Bbbk [z_{124134}] =  \toba_5^{\co \pi_2}$, so $X'_5:= \mH^{\co \pi_2\# \id} =\Bbbk [z_{124134}]\# 1$; set $X_5=\Ss(X'_5)$. There exists an algebra map $\phi \colon X_5 \to \cA_5$, $\phi(z_{124134}g_1^{-12}g_4^{-6})=\gamma_4(z_{124134})g_1^{-12}g_4^{-6}-\lambda_4g_1^{-12}g_4^{-6}$, which is $\mH_5$-colinear. We claim that this map is also $\mH_5$-linear. Indeed, 
$\phi(g_i \cdot z_{124134}g_1^{-12}g_4^{-6}) = g_i \cdot \phi(z_{124134}g_1^{-12}g_4^{-6})$ by \eqref{eq:lambda-conditions}. 
Set $\ytt_{124134}=\gamma_4(z_{124134})$.
By Lemma \ref{lem:Z-HV-str} \ref{item:Z-HV-central},
$\phi(x_i \cdot z_{124134}g_1^{-12}g_4^{-6})=0$. Notice that 
\begin{align*}
x_i \cdot \phi \big(z_{124134}g_1^{-12}g_4^{-6}\big) &= \big(y_i \ytt_{124134} - \chi_1^{12}\chi_4^{6}(g_i) \ytt_{124134}y_i\big)g_1^{-12}g_4^{-6}.
\end{align*}
By direct computation, 
\begin{align*}
\rho\big(y_i \ytt_{124134} - \chi_1^{12}\chi_4^{6}(g_i) \ytt_{124134}y_i\big) &= \big(y_i \ytt_{124134} - \chi_1^{12}\chi_4^{6}(g_i) \ytt_{124134}y_i\big)\ot 1,
\end{align*}
hence $y_i \ytt_{124134} - \chi_1^{12}\chi_4^{6}(g_i) \ytt_{124134}y_i =\mu$ for some $\mu_i\in\Bbbk$, $i \in \I_4$.

On the other hand, $\gamma_4$ is $\Gamma$-linear by \cite[Proposition 5.8 (c)]{A+}. If $i \in \I_3$, 
\begin{align*}
g_i \cdot\big(y_i \ytt_{124134} - \chi_1^{12}\chi_4^{6}(g_i) \ytt_{124134}y_i\big)&=
-q_1^2\big(y_i \ytt_{124134} - \chi_1^{12}\chi_4^{6}(g_i) \ytt_{124134}y_i\big),
\\
g_4\cdot\big(y_i \ytt_{124134} - \chi_1^{12}\chi_4^{6}(g_i) \ytt_{124134}y_i\big)&=
q_2^5\omega \big(y_i \ytt_{124134} - \chi_1^{12}\chi_4^{6}(g_i) \ytt_{124134}y_i\big).
\end{align*}
If $\mu_i \neq 0$ we get $-q_1^2=1=q_2^5\omega$, thus $\omega=\omega^{10}=(q_1q_2)^{10}=-\omega$, a contradiction. Hence $\mu_i=0$ for $i \in \I_3$. Analogously, if $\mu_4\neq 0$, then $-q_2^4=1=q_1^3$, a contradiction since $q_1q_2\in \Gp_6$. Thus $\phi$ is $\mH_5$-linear.

By \cite[Theorem 4]{Gu} $\cA$ is a $\mH$-cleft object. The claim about the section $\gamma$ follows from \cite[Proposition 5.8]{A+}.
\epf

\subsection{Computing the liftings} 

Let $\bsl \in \cR_{\nombre}$ and define $\cL_1=\mH_1$, but we change the labels of the generators to $(a_i)_{i\in\I_4}$,
\begin{align*}
\cL_2(\bsl) & = \cL_1(\bsl) / \left\langle a_1^2 - \lambda_1 ( 1 - g_1^2)\right\rangle,
\\
\cL_3(\bsl) & = \cL_2(\bsl) / \left\langle  a_1 a_2 + a_3 a_1 + a_2 a_3 - \lambda_2 (1 - g_1 g_2 )\right\rangle,
\\
\cL_4(\bsl) & = \cL_3(\bsl) / \left\langle  a_{214} - \omega a_{134} - \lambda_2 a_4 \right\rangle.
\end{align*}

Notice that, in $\cL_2(\bsl)$,
\begin{align*}
a_i^2 - \lambda_1 ( 1 - g_i^2) &= g_{5-i} (a_1^2 - \lambda_1 ( 1 - g_1^2))g_{5-i}^{-1} =0, & &i\in\I_{2,3}.
\end{align*}
The same happens for the other sets of relations in $\cL_3(\bsl)$ and $\cL_4(\bsl)$.

\begin{lemma}\label{lem:liftingsHV} Let $\bsl \in \cR_{\nombre}$, $i\in\I_4$.
Then $\cA_i (\bsl)$ is a $(\cL_i(\bsl),\mH_i)$-biGalois object.
\end{lemma}

\pf
For $i=1,2,3$, $\cL_i(\bsl) \simeq L(\cA_i(\bsl),\mH_i)$ by \cite[Corollary 5.12]{A+}, since all the involved relations are skew-primitive elements of $\mT$. 

For $i=4$, we use Remark \ref{rem:comultiplication-relations} to compute
\begin{align*}
\rho_3\big( y_{214}-\omega y_{134} - \lambda_2 y_4 \big) & = \big( a_{214}-\omega a_{134} - \lambda_2 a_4 \big) \ot 1 
\\ & + g_2g_1g_4 \ot \big( y_{214}-\omega y_{134} - \lambda_2 y_4 \big).
\end{align*}
Hence $\cL_4(\bsl) \simeq L(\cA_4(\bsl),\mH_4)$, again by \cite[Corollary 5.12]{A+}.
\epf

Thanks to the previous Lemma and \cite[Corollary 5.12]{A+}, there exists a $(g_1^{12}g_4^6,1)$-primitive element $\att_{124134}\in\cL_4$
\footnote{We performed an algorithm in GAP to compute explicitly $\att_{124134}$ and $\gamma_4(z_{124134})$ but the program did not finish the computation.} such that
\begin{align}\label{eq:a124134-defn}
\att_{124134}\ot 1 = \delta_4\big(\gamma_4(z_{124134})\big)-g_1^{12}g_4^6 \ot \gamma_4(z_{124134}).
\end{align}
Next we define
\begin{align}\label{eq:L-lambda}
\cL(\bsl) := \cL_4(\bsl)/ \langle a_4^6- \lambda_3(1-g_4^6), \att_{124134} - \lambda_4 (1-g_1^{12}g_4^6) \rangle.
\end{align}

\begin{theorem}\label{thm:liftingsHV} 
Let $\Gamma$ be a group with a principal realization of $\nombre$ as in \S \ref{subsec:Gamma3}. Let $\bsl \in \cR_{\nombre}$, see \eqref{eq:lambda-conditions}. Then
\begin{enumerate}[leftmargin=*,label=\rm{(\alph*)}]
\item\label{item:liftings-bigalois} $\cL (\bsl)\simeq L(\cA(\bsl),\toba(\nombre)\#\Bbbk\Gamma)$.
\item\label{item:liftings-lifting} $\cL (\bsl)$ is a lifting of $\toba(\nombre)$ over $\Bbbk\Gamma$.
\item\label{item:liftings-cocycle} $\cL(\bsl)$ is a cocycle deformation of $\toba(\nombre)\#\Bbbk\Gamma$.
\end{enumerate}

Conversely, if $L$ is lifting of $\toba(\nombre)$ over $\Bbbk\Gamma$, then there exist $\bsl \in \cR_{\nombre}$ 
such that $L \simeq \cL (\bsl)$.
\end{theorem}
\pf
First, \ref{item:liftings-bigalois} follows again by \cite[Corollary 5.12]{A+} since $z_4$ and $z_{124134}$ are skew primitive in $\mH_4$, and this implies that \ref{item:liftings-cocycle} holds. Now \ref{item:liftings-lifting} is a consequence of Lemma \ref{lem:A-cleft} and \cite[Proposition 4.14. (c)]{A+}.

\medspace

Conversely, assume that $L$ is a lifting of $\toba(\nombre)$ over $\Bbbk\Gamma$. Let $\phi: \mT= T(\nombre)\#\Bbbk\Gamma \twoheadrightarrow L$ be a lifting map as in \cite[Proposition 2.4]{AV}. As in \cite[Theorem 3.5]{AAG}, we shall attach to $\phi$ a family of scalars $\bsl \in \cR_{\nombre}$ as follows. Since $\nombre \# \Gamma \subset \sum_{i\in\I_4, \, g \in \Gamma} \ku g_ig \wedge \ku g$, it follows that the first term of the coradical filtration of $L$ satisfies
\begin{align}\label{eq:lifting-coradical-1}
L_1 = \phi (\ku\Gamma \oplus \nombre\#\ku\Gamma) = \ku\Gamma + \sum_{i\in\I_4, \, g \in \Gamma} \ku g_ig \wedge \ku g.
\end{align}

Let $r \in \cG_0\cup\cG_1\cup\cG_2$. Since $r \in \mT$ is $(g_r, 1)$-primitive,  \eqref{eq:lifting-coradical-1} implies that $\phi(r) \in \Bbbk(1-g_r) + \ku g_r\wedge \ku \subset L_1$. By Lemma \ref{lem:gr-neq-gi}, $g_r\neq g_i$ for all $i\in\I_4$.
Thus $\phi(r)=\lambda_r(1-g_r)$ for some $\lambda_r\in\Bbbk$. 

Let $r=a_4 a_{h4} - q_2 a_{h4} a_4\in\cG_0$, for $h\in\I_3$. Suppose that $\lambda_r\neq 0$. Then 
\begin{align*}
1&=\chi_h\chi_4^2(g_4)=\omega q_2, & 1&=\chi_h\chi_4^2(g_h)=-q_1^2,
\end{align*}
which implies that $q_1q_2\in\Gp_{12}$, a contradiction. Hence $\lambda_r=0$.

Next we consider $a_h^2\in\cG_1$: we denote $\phi(a_h^2)=\lambda_{1h}(1-g_h^2)$. For $h\in\I_{2,3}$,
\begin{align*}
\lambda_{11}(1-g_1^2) & = \phi(a_1^2) = \phi(g_{5-h} a_h g_{5-h}^{-1}) = g_{5-h}\phi(a_h) g_{5-h}^{-1} = \lambda_{1h}(1-g_h^2),
\end{align*}
where we use that $g_1^2=g_h^2$. Hence $\lambda_{11}=\lambda_{1h}$, and we denote the common scalar simply by $\lambda_{1}$. A similar computation shows that $\lambda_{1}$ satisfies \eqref{eq:lambda-conditions}.

A similar argument shows that $\phi(a_1 a_2 + a_3 a_1 + a_2 a_3) = \lambda_2 (1 - g_1 g_2)$ for some $\lambda_2\in\Bbbk$ satisfying \eqref{eq:lambda-conditions}. Thus $\phi$ descends to a Hopf algebra map $\phi:\cL_3(\bsl)\twoheadrightarrow L$ for any choice of $\lambda_3$ and $\lambda_4$.

Next, $a_{214} - \omega a_{134} - \lambda_2 a_4\in \cL_3(\bsl)$ is $(g_1g_3g_4,1)$-primitive. By Lemma \ref{lem:gr-neq-gi} there exists $\mu\in\Bbbk$ such that $\phi(a_{214} - \omega a_{134} - \lambda_2 a_4)=\mu(1-g_1g_3g_4)$. Suppose that $\mu\neq 0$. Conjugation by $g_1$ and $g_4$ (as we applied for the previous relations) imply that $q_1=1$ and $-q_2^2\omega^2=1$ respectively, but this is a contradiction since $q_1q_2=-\omega$. Thus $\mu=0$ and $\phi$ descends to a Hopf algebra map $\phi:\cL_4(\bsl)\twoheadrightarrow L$ for any choice of $\lambda_3$ and $\lambda_4$.

Finally, $a_4^6$ is $(g_4^6,1)$-primitive and $\att_{124134}$ is $(g_1^{12}g_4^6,1)$-primitive. An argument analogous to the previous relations shows that 
\begin{align*}
\phi(a_4^6)&=\lambda_3(1-g_4^6), &
\phi(\att_{124134})&=\lambda_4(1-g_1^{12}g_4^6),
\end{align*}
for some $\lambda_3, \lambda_4\in\Bbbk$ satisfying \eqref{eq:lambda-conditions}. Thus $\phi$ descends to a Hopf algebra map $\phi:\cL(\bsl)\twoheadrightarrow L$; as the restriction of $\phi$ to the first term $\cL(\bsl)_1$ of the coradical filtration is injective, $\phi$ is an isomorphism by \cite[Theorem 5.3.1]{Mo}. 
\epf

\end{document}